\documentclass[11pt,a4paper,reqno]{amsart}
\usepackage{calc,%xcolor,
amsfonts,amsthm,amscd,epsfig,
psfrag,amsmath,amssymb,enumerate,
graphicx}

\reversemarginpar
\newlength\fullwidth
\numberwithin{equation}{section}

\DeclareMathSymbol{\leqslant}{\mathalpha}{AMSa}{"36} % nicer `smaller or equal'
\DeclareMathSymbol{\geqslant}{\mathalpha}{AMSa}{"3E} % nicer `larger or equal'
\DeclareMathSymbol{\eset}{\mathalpha}{AMSb}{"3F}     % nicer `emptyset'
\renewcommand{\leq}{\;\leqslant\;}                   % redef. of < or =
\renewcommand{\geq}{\;\geqslant\;}                   % redef. of > or =
\newcommand{\eps}{\epsilon}

\def\1{\ifmmode {1\hskip -3pt \rm{I}} \else {\hbox {$1\hskip -3pt \rm{I}$}}\fi}

\newtheorem{Theorem}{Theorem}[section]
\newtheorem{Lemma}[Theorem]{Lemma}
\newtheorem{Proposition}[Theorem]{Proposition}
\newtheorem{Corollary}[Theorem]{Corollary}
\newtheorem{remark}[Theorem]{Remark}

\newtheorem{claim}[Theorem]{Claim}
\newtheorem{definition}[Theorem]{Definition}

\newcommand{\cB}{\ensuremath{\mathcal B}}
\newcommand{\cC}{\ensuremath{\mathcal C}}
\newcommand{\cD}{\ensuremath{\mathcal D}}

\newcommand{\cF}{\ensuremath{\mathcal F}}

\newcommand{\cK}{\ensuremath{\mathcal K}}
\newcommand{\cL}{\ensuremath{\mathcal L}}

\newcommand{\cN}{\ensuremath{\mathcal N}}

\newcommand{\bbE}{{\ensuremath{\mathbb E}} }
\newcommand{\bbF}{{\ensuremath{\mathbb F}} }

\newcommand{\bbN}{{\ensuremath{\mathbb N}} }

\newcommand{\bbP}{{\ensuremath{\mathbb P}} }

\newcommand{\bbZ}{{\ensuremath{\mathbb Z}} }

\newcommand{\ent}{{\rm Ent} }

\let\a=\alpha \let\b=\beta   \let\d=\delta  \let\e=\varepsilon
 \let\g=\gamma \let\h=\eta      \let\l=\lambda
      \let\o=\omega      
  \let\s=\sigma \let\t=\tau   
  
   \let\G=\Gamma  \let\L=\Lambda 
\let\O=\Omega      

\def\smallno{\smallskip\noindent}
\def\medno{\medskip\noindent}
\def\bigno{\bigskip\noindent}
\def\\{\hfill\break}

\def\thsp{\thinspace}

\def\tthsp{\kern .083333 em}

\def\?{\mskip -10mu}
\def\indbox#1{\hbox to \parindent{\hfil\ #1\hfil} }

\newfam\msafam
\newfam\msbfam
\newfam\eufmfam
\def\hexnumber#1{%
  \ifcase#1 0\or 1\or 2\or 3\or 4\or 5\or 6\or 7\or 8\or
  9\or A\or B\or C\or D\or E\or F\fi}
\font\tenmsa=msam10
\font\sevenmsa=msam7
\font\fivemsa=msam5
\textfont\msafam=\tenmsa
\scriptfont\msafam=\sevenmsa
\scriptscriptfont\msafam=\fivemsa
\edef\msafamhexnumber{\hexnumber\msafam}%
\mathchardef\restriction"1\msafamhexnumber16
\mathchardef\ssim"0218
\mathchardef\square"0\msafamhexnumber03
\mathchardef\eqd"3\msafamhexnumber2C
\def\QED{\ifhmode\unskip\nobreak\fi\quad
  \ifmmode\square\else$\square$\fi}
\font\tenmsb=msbm10
\font\sevenmsb=msbm7
\font\fivemsb=msbm5
\textfont\msbfam=\tenmsb
\scriptfont\msbfam=\sevenmsb
\scriptscriptfont\msbfam=\fivemsb

\font\teneufm=eufm10
\font\seveneufm=eufm7
\font\fiveeufm=eufm5
\textfont\eufmfam=\teneufm
\scriptfont\eufmfam=\seveneufm
\scriptscriptfont\eufmfam=\fiveeufm

\def\({\left(}
\def\){\right)}
\let\neper=e
\let\ii=i

\def\ie{\hbox{\it i.e.\ }}
\let\id=\identity
\let\emp=\emptyset
\let\sset=\subset

\let\setm=\backslash
\def\nep#1{ \neper^{#1}}

\let\imp=\Rightarrow

\def\tc{\thsp | \thsp}
\let\<=\langle
\let\>=\rangle

\def\Var{ \mathop{\rm Var}\nolimits }

\def\gap{\mathop{\rm gap}\nolimits}

\def\scalprod#1#2{ \thsp<#1, \thsp #2>\thsp }
\def\inte#1{\lfloor #1 \rfloor}

\outer\def\nproclaim#1 [#2]#3. #4\par{\medbreak \noindent
   \talato(#2){\bf #1 \Thm[#2]#3.\enspace }%
   {\sl #4\par }\ifdim \lastskip <\medskipamount
   \removelastskip \penalty 55\medskip \fi}
\def\thmm[#1]{#1}
\def\teo[#1]{#1}
\def\sttilde#1{%
\dimen2=\fontdimen5\textfont0
\setbox0=\hbox{$\mathchar"7E$}
\setbox1=\hbox{$\scriptstyle #1$}
\dimen0=\wd0
\dimen1=\wd1
\advance\dimen1 by -\dimen0
\divide\dimen1 by 2
\vbox{\offinterlineskip%
   \moveright\dimen1 \box0 \kern - \dimen2\box1}
}
\begin{document}

\title[Kinetically constrained models]{Kinetically constrained spin models}

\author[N. Cancrini]{N. Cancrini}
\email{nicoletta.cancrini@roma1.infn.it}
\address{Dip. Matematica Univ.l'Aquila, 1-67100 L'Aquila, Italy}

\author[F. Martinelli]{F. Martinelli}
\email{martin@mat.uniroma3.it}
\address{Dip.Matematica, Univ. Roma Tre, Largo S.L.Murialdo 00146, Roma, Italy}

\author[C. Roberto]{C. Roberto}
\email{cyril.roberto@univ-mlv.fr}
\address{L.A.M.A., Univ Marne-la-Vall\'ee, 5 bd Descartes 77454 Marne-la-Vall\'ee France}
\author[C. Toninelli]{C. Toninelli}
\email{Cristina.Toninelli@lpt.ens.fr}
\address{Laboratoire de Probabilit\'es et Mod\`eles Al\`eatoires
  CNRS-UMR 7599 Universit\'es Paris VI-VII 4, Place Jussieu F-75252 Paris Cedex 05 France}

\thanks{We would like to thank A.~Gandolfi and J.~van~den Berg for
  a useful discussion on high dimensional percolation and H.~C. Andersen
  for some enlightening correspondence}

\begin{abstract}
  We analyze the density and size dependence of the relaxation time for
  kinetically constrained spin models (KCSM) intensively studied in the
  physical literature as simple models sharing some of the features of a
  glass transition.  KCSM are interacting particle systems on $\bbZ^d$
  with Glauber-like dynamics, reversible w.r.t. a simple product i.i.d
  Bernoulli($p$) measure.  The essential feature of a KCSM is that the
  creation/destruction of a particle at a given site can occur only if
  the current configuration of empty sites around it satisfies certain
  constraints which
  completely define each specific model. No other interaction is present
  in the model. From the mathematical point of view, the basic issues
  concerning positivity of the spectral gap inside the ergodicity region
  and its scaling with the particle density $p$ remained open for most
  KCSM (with the notably exception of the East model in $d=1$
  \cite{Aldous}). Here for the first time we: i) identify the ergodicity
  region by establishing a connection with an associated bootstrap
  percolation model; ii) develop a novel multi-scale approach which
  proves positivity of the spectral gap in the whole ergodic region;
  iii) establish, sometimes optimal, bounds on the behavior of the
  spectral gap near the boundary of the ergodicity region and iv)
  establish pure exponential decay for the persistence function (see
  below).  Our techniques are flexible enough to allow a variety of
  constraints and our findings disprove certain conjectures which
  appeared in the physical literature on the basis of numerical
  simulations.

\bigno

\end{abstract}

\bigno

%\footnote{{????}}

%\date{22 Marzo 06}
\maketitle

\thispagestyle{empty}

\section{Introduction}
Kinetically constrained spin models (KCSM) are interacting particle
systems on the integer lattice $\bbZ^d$. A configuration is defined by
assigning to each site $x$ its occupation variable $\eta_x\in\{0,1\}$.
The evolution is given by a simple Markovian stochastic dynamics of
Glauber type with generator $\cL$. Each site waits an independent, mean
one, exponential time and then, provided that the current configuration
around it satisfies an apriori specified constraint which does not involve
$\h_x$, it refreshes its state by declaring it to be occupied with
probability $p$ and empty with probability $q=1-p$. Detailed balance
w.r.t. Bernoulli($p$) product measure $\mu$ is easily verified and $\mu$
is therefore an invariant reversible measure for the process.

These models have been introduced in physical literature \cite{FA1,FA2}
to model liquid/glass transition and more generally the slow ``glassy''
dynamics which occurs in different systems (see \cite{Ritort,TB} for
recent review). In particular, they were devised to mimic the fact that
the motion of a molecule in a dense liquid can be inhibited by the
presence of too many surrounding molecules.  That explains why, in all
physical models, the constraints specify the maximal number of particles
on certain sites around a given one in order to allow
creation/destruction on the latter.  As a consequence, the dynamics of
KCSM becomes increasingly slow as $p$ is increased.  Moreover there
usually exist configurations with all creation/destruction rates identically equal to zero
(blocked configurations), a fact that implies the existence
of several invariant measures (see \cite{Lalley} for a somewhat detailed
discussion of this issue in the context of the North-East
model) and produce unusually long
mixing times compared to standard high-temperature stochastic Ising models (see
section \ref{entropy} below). Finally we observe that a KCSM model is in
general not attractive so that the usual coupling arguments valid for
e.g. ferromagnetic stochastic Ising models cannot be applied.

The above  little discussion explains why the basic issues concerning
the large time behavior of the process, even if started from the
equilibrium reversible measure $\mu$, are non trivial and justifies why
they remained open for most of the interesting models, with the only exception
of the East model \cite{Aldous}. This is a one-dimensional model for which
creation/destruction at a given site can occur only if the nearest
neighbor to its right is empty. In \cite{Aldous} it has been
proved that the generator $\cL$ of the East model has a positive
spectral gap for all $q>0$, which, for $q\downarrow 0$, shrinks faster than any
polynomial in $q$ (see section \ref{specificmodels} for more details). However, the
method in \cite{Aldous} uses quite heavily
the specifics of the model and
its extension to higher dimensions or to other models introduced in
physical literature seems to be non trivial.  Among the latter we just
recall the North-East model (N-E)  \cite{RJM} in $\bbZ^2$ and the
Fredrickson Andersen $j\le d$ spin facilitated (FA-jf) \cite{FA1} models
in $\bbZ^d$. For the first, destruction/creation at a given site can
occur only if its North and East neighbors are empty, while for the
FA-jf model the constraint requires that at least $j$ among the nearest
neighbors are empty.

The main achievements of this paper can be described as follows.  In
section 2.3, given a generic KCSM with constraints satisfying few rather
mild conditions, we first identify the critical value of the density of
vacancies $q_c=\inf\{q: 0 \text{ is a simple eigenvalue of $\cL$}\}$
with the critical value of a naturally related bootstrap percolation
model. Notice that a general result on Markov semigroups (see Theorem
\ref{ergodic} below) implies that for any $q>q_c$ the reversible measure
$\mu$ is mixing for the process generated by $\cL$. Next, in section 3,
we identify a natural general condition on the associated bootstrap
percolation model which implies the positivity of the spectral gap of
$\cL$. In its simplest form the condition requires that the probability
that a large cube is internally spanned (\ie the block does not contain
blocked configurations, see definition \ref{IS} below) is close to one.  For all the models
discussed in section \ref{specificmodels} our condition is satisfied for
all $q$ strictly larger than $q_c$. Our findings disprove some
conjectures appeared in the physical literature \cite{GPG,Ha}, based on
numerical simulations and approximated analytical treatments, on the
existence of a second critical point $q_c'>q_c$ at which the spectral
gap vanishes. The main ingredients in the proof are multi-scale
arguments, the bisection technique of \cite{SFlour} combined with the
novel idea of considering auxiliary constrained models on large length
scales with scale dependent constraints (see sections \ref{GM} and \ref{Proof of theorem 1}) . At the end of the section we
also analyze the so called persistence function $F(t)$ which represents
the probability for the equilibrium process that the occupation variable
at the origin does not change before time $t$. We prove that, whenever
the spectral gap is strictly positive, $F(t)$ must decay
exponentially. This, together with the above results, disproves previous
conjectures of a stretched exponential decay of the form $F(t)\simeq
\exp(-t/\tau)^{\beta}$ with $\beta<1$ for FA1f in $d=1$ \cite{BG,WBG1}
and for FA2f in $d=2$ \cite{Ha} \footnote{For a different Ising-type
  constrained model in which the kinetic constraint prevents spin-flip
  which do not conserve the energy, Spohn \cite{Spohn} has proved long
  ago that the time autocorrelation of the spin at the origin decays as
  a stretched exponential}. For the North-East model at the critical
point we show instead (see corollary \ref{slow}) that $\int_0^\infty dt
\, F(\sqrt{t})=\infty$, a signature of a slow polynomial decay.

After establishing the positivity of the spectral gap, in section
\ref{specificmodels} we analyze its behavior as $q\downarrow q_c$ for
some of the models discussed in section 2.3.  For the East model
($q_c=0$) we significantly improve the lower bound on the spectral gap
proved in \cite{Aldous} and claimed to provide the leading behavior in
\cite{SE}. Our lower bound, in leading order, coincides with the upper
bound of \cite{Aldous}, yielding that the gap shrinks as
$q^{\log_2(q)/2}$ for small values of $q$.

For the FA-1f model ($q_c=0$) we show that for $q\approx 0$, the
spectral gap is $O(q^3)$
in $d=1$, $O(q^2)$ in $d=2$ apart from logarithmic corrections
and between $O(q^{1+2/d})$ and $O(q^2)$ in $d\ge 3$. Again these findings
disprove previous claims in $d=2,3$ \cite{WBG1} .

For the FA-2f model ($q_c=0$) in e.g. $d=2$ we get
instead
\begin{equation}
  \label{eq:0.1}
\exp(-c/q^5)\le \gap(\cL) \le
\exp\Bigl(-\frac{\pi^2}{18q}\bigl(1+o(1)\bigr)\Bigr)
\end{equation}
as $q\downarrow 0$.  Notice that the r.h.s. of \eqref{eq:0.1}
represents the inverse of the critical length for bootstrap percolation
\cite{Holroyd}, \ie the smallest length scale above which a region of
the lattice becomes mobile or unjammed under the FA-2f dynamics, and it
has been conjectured \cite{R,TBF2} to provide the
leading behavior of the spectral gap for small values of $q$.

As explained above, the techniques developed in this paper are flexible
enough to deal with a variety of KCSM even on more general graphs
\cite{preparation0} and, possibly, with some non trivial interaction
between the occupation variables. Furthermore it seems that they could
also be applied to kinetically constrained models with Kawasaki
(i.e. conservative) rather than Glauber dynamics.

\section{The general models}
\label{general models}
\subsection{Setting and notation}
\label{setting}
The models considered here are defined on the integer lattice
$\bbZ^d$ with sites $x = (x_1,\dots,x_d)$ and basis vectors $\vec e_1=(1,\dots,0),
\vec e_2=(0,1,\dots,0),\dots, \vec e_d=(0,\dots,1)$.
On $\bbZ^d$ we will consider the Euclidean norm $\|x\|$, the $\ell^1$
(or graph theoretic) norm $\|x\|_1$ and the sup-norm $\|x\|_\infty$. The
associated distances will be denoted by
$d(\cdot, \cdot)$, $d_1(\cdot, \cdot)$ and $d_\infty(\cdot, \cdot)$ respectively.
For any vertex $x$ we let
$$
\begin{array}{l}
  \cN_x=\{y\in \bbZ^d:\ d_1(x,y)=1\}, \\
  \cK_x=\{y\in \cN_x:\ y=x+\sum_{i=1}^d \a_i\vec e_i,\ \a_i\ge 0\}\\
  \cN_x^*=\{y\in \bbZ^d:\ y=x+\sum_{i=1}^d \a_i\vec e_i,\
  \a_i=\pm 1,0 \text{ and }\sum_i\a^2_i\neq 0\} \\
  \cK^*_x=\{y\in \cN_x^*:\ y=x+\sum_{i=1}^d \a_i\vec e_i,\ \a_i=1,0\}
\end{array}
$$
and write $x\sim y$ if $y\in \cN_x^*$.
\begin{figure}[h]
\psfrag{x}{$x$}
\psfrag{N}{${\mathcal N}_x$}
\psfrag{K}{${\mathcal K}_x$}
\psfrag{O}{${\mathcal N}_x^*$}
\psfrag{L}{${\mathcal K}_x^*$}
\includegraphics[width=.80\columnwidth]{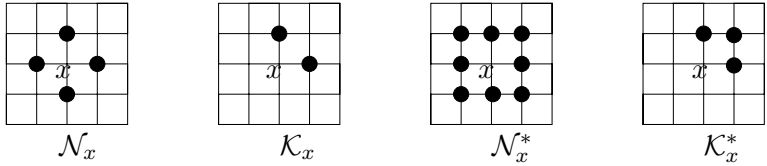}
\caption{The various neighborhoods of a vertex $x$ in two dimensions}
\end{figure}
The neighborhood, the *-neighborhood, the oriented and *-oriented neighborhoods $\partial
\L,\ \partial^*\L, \partial_+ \L,\ \partial_+^* \L$ of a finite subset $\L\sset
\bbZ^d$ are defined accordingly as
$
\partial \L:=\{\cup_{x\in \L}\cN_x\}\setminus \L,\
\partial^* \L:=\{\cup_{x\in \L}\cN^*_x\}\setminus \L,\
\partial_+\L:=\left\{\cup_{x\in \L}\cK_x\right\}\setminus \L,\
\partial_+^* \L:=\left\{\cup_{x\in \L}\cK^*_x\right\}\setminus \L.
$
A rectangle $R$ will be a set of sites of the form
$$
R:=[a_1,b_1]\times\dots\times [a_d,b_d]
$$
while the collection of finite subsets of $\bbZ^d$ will be denoted by
$\bbF$.

\medno
The pair $\left(S,\nu\right)$ will denote a finite probability space with
  $\nu(s)>0$ for any $s\in S$.
$G\sset S$ will denote a distinguished event in $S$, often referred to as
the set of ``good states'', and
$q\equiv \nu(G)$ its probability.

Given $\left(S,\nu\right)$ we will consider the configuration space
$\O=S^{\bbZ^d}$ equipped with the product measure $\mu:=\prod_{x\in
  \bbZ^d}\nu_x$, $\nu_x\equiv \nu$. Similarly we define $\O_\L$ and
$\mu_\L$ for any subset $\L\sset \bbZ^d$. Elements of $\O$ ($\O_\L$) will be denoted
by Greek letters $\o,\h$ ($\o_\L,\h_\L$) etc while the variance w.r.t
$\mu$
by $\Var$ ($\Var_\L$). Finally we will use the shorthand notation
$\mu(f)$ to denote the expected value of any $f\in L^1(\mu)$.
%The set of functions $f:\O\mapsto \bbR$ which
%depend on finitely many variables will be denoted by $\cD$.

\subsection{The Markov process}
\label{Markov preocess}
The interacting particle models that will be studied here are
Glauber type Markov processes in $\O$, reversible w.r.t. the measure
$\mu$ (or $\mu_\L$ if considered in $\O_\L$) and characterized by a
collection of
\emph{influence classes} $\{\cC_x\}_{x\in \bbZ^d}$ and by a choice of the good event $G\subset S$.
For any $x$, $\cC_x$ is a collection of subsets of $\bbZ^d$ (see below for some of the
most relevant examples).
The collection of influence classes will satisfy the following basic hypothesis:
\begin{enumerate}[a)]
\item \emph{independence of $x$:} for all $x\in \bbZ^d$ and all $A\in
  \cC_x$ $x\notin A$ ;
\item \emph{translation invariance:} $\cC_x=\cC_0+x$ for all $x$;
\item \emph{finite range interaction:} there exists $r<\infty$ such that
  any element of
  $\cC_x$ is contained in $\cup_{j=1}^r\{y:\ d_1(x,y)=j\}$
\end{enumerate}
\begin{definition}
Given a vertex $x\in \bbZ^d$ we will say that the
constraint at $x$ is satisfied by
 the configuration $\o$ if the indicator
\begin{equation*}
  c_{x}(\o)=
  \begin{cases}
 1 & \text{if there exists a set $A\in
\cC_x$ such that $\o_y\in G$ for all $y\in A$}\\
0 & \text{otherwise}
  \end{cases}
\end{equation*}
is equal to one.
\end{definition}
The process that will be studied in the sequel can be informally described
as follows. Each vertex $x$ waits an independent mean one exponential
time and then, provided that the current configuration $\o$ satisfies
the constraint at $x$, the value $\o_x$ is refreshed with a
new value in $S$ sampled from $\nu$ and the all procedure starts again.

The generator $\cL$ of the process can be constructed in a standard way
(see e.g. \cite{Liggett,Lalley}) and it is a non negative self-adjoint
operator on $L^2(\O,\mu)$ with domain $Dom(\cL)$ and Dirichlet form given by
$$
\cD(f)=\sum_{x\in \bbZ^d}\mu\(c_{x} \Var_x(f)\),\quad f\in Dom(\cL)
$$
Here $\Var_{x}(f)\equiv\int d\nu(\o_x) f^2(\o)
- \(\int d\nu(\o_x)f(\o)\)^2$ denotes the local variance with respect to
the variable $\o_x$ computed while the other variables are held fixed.
To the generator $\cL$ we can associate the Markov semigroup
$P_t:=\nep{t\cL}$ with reversible invariant measure $\mu$.

Notice that the constraints $c_x(\o)$ are increasing functions w.r.t the
partial order in $\O$ for which $\o\le\o'$ iff $\o'_x\in G$ whenever
$\o_x\in G$.
However that does not imply in general that the process generated by
$\cL$ is attractive in the sense of Liggett \cite{Liggett}.

Due to the fact that the jump rates are not bounded away from zero, the
reversible measure $\mu$ is certainly not the only invariant measure
(there exists initial configurations that are blocked forever) and an
interesting question is therefore whether $\mu$ is ergodic or mixing for
the Markov process and whether there exist other translation invariant,
ergodic stationary measures. To this purpose it is useful to recall the following
well known result (see e.g. Theorem 4.13 in \cite{Liggett}).
\begin{Theorem}
\label{ergodic}
  The following are equivalent,
  \begin{enumerate}[(a)]
  \item $\lim_{t\to \infty} P_tf=\mu(f)$ in $L^2(\mu)$ for all
    $f\in L^2(\mu)$.
\item $0$ is a simple eigenvalue for $\cL$.
  \end{enumerate}
\end{Theorem}
Clearly $(a)$ above implies that $\lim_{t\to
  \infty}\mu\left(fP_tg\right)=\mu(f)\mu(g)$ for any $f,g\in L^2(\mu)$,
\ie $\mu$ is mixing and therefore ergodic.
\begin{remark}
Even if $\mu$ is mixing there will exist in
general infinitely many stationary measures, \ie probability measures
$\tilde \mu$ satisfying $\tilde\mu P_t=\tilde\mu$ for all $t\ge 0$. As
an example take an arbitrary probability measure $\tilde \mu$ such that $\tilde\mu\bigl(\{S\setminus
G\}^{\bbZ^d}\bigr)=1$.  We refer the interested reader to \cite{Lalley}
for a discussion of this point in the context of the North-East model
(see below).
\end{remark}

In a finite region $\L\sset \bbZ^d$ the process, a continuous time
Markov chain in this case, can be defined analogously but some care has
to be put in order to correctly define the constraints $c_x$ for those $x\in
\L$ such that their influence class $\cC_x$ is not entirely contained
inside $\L$.

One possibility is to modify in a $\L$-dependent way the definition of
the influence classes e.g. by defining
$$
\cC_{x,\L}:=\{A\cap \L;\ A\in \cC_x\}
$$
Although such an approach is feasible and natural, at least for some of the models
discussed below, an important drawback is a loss of ergodicity of the
chain. One is then forced to consider the chain restricted to an ergodic
component making the whole analysis more cumbersome (see section \ref{Some further observations}).

Another alternative is to imagine that the
configuration $\o$ outside $\L$ is frozen and equal to some reference
configuration $\t$ that will be referred to as the \emph{boundary
  condition} and to define the finite volume constraints with boundary
condition $\t$ as
$$
c^\t_{x,\L}(\o_\L):= c_x(\o_\L\cdot\t),
$$
where $\o_\L\cdot \t$ simply denotes the configuration equal
to $\o_\L$ inside $\L$ and equal to $\t$ outside. Since we want the
Markov chain to be ergodic $\t$ will need to be in the good set $G$ for some of the
vertices outside $\L$.  Instead of discussing this issue in a very
general context we will now describe the basic models and solve the
problem of boundary conditions for each one of them.

\subsection{0-1 Kinetically constrained spin models}
\label{models}
In most models considered in the physical literature the finite
probability space $(S,\nu)$ is simply the two state-space $\{0,1\}$
and the good set $G$ is conventionally chosen as the empty state
$\{0\}$. Any model with these features will be called a ``0-1 KCSM''
(kinetically constrained spin model).

Given a 0-1 KCSM, the parameter $q=\mu(\h_0=0)$ can be varied in $[0,1]$ while
keeping fixed the
basic structure of the model (\ie the notion of the good set and the
functions $c_x$'s expressing the constraints) and it is natural to define a critical
value $q_c$ as
$$
q_c=\inf\{q\in[0,1]:\, 0 \text{ is a simple eigenvalue of $\cL$}\}
$$
As we will prove below $q_c$ coincides with the \emph{bootstrap
  percolation threshold} $q_{bp}$
of the model defined as follows \cite{Schonmann} \footnote{In most of
  the boostrap percolation literature the role of the $0$'s and the
  $1$'s is inverted}. For any $\h\in \O$
define the bootstrap map $T:\O\mapsto \O$ as
\begin{equation}
  \label{eq:bootstrap map}
  T(\h)_x=0 \quad \text{if either}\quad
\h_x=0 \quad \text{or}\quad c_x(\h)=1.
\end{equation}
Denote by $\mu^{(n)}$ the probability
measure on $\O$ obtained by iterating $n$-times the above mapping
starting from $\mu$. As $n\to \infty$ $\mu^{(n)}$ converges to
a limiting measure $\mu^{(\infty)}$ \cite{Schonmann} and it is natural
to define the critical value $q_{bp}$ as
$$
q_{bp}=\inf\{q\in [0,1]:\,\mu^{(\infty)}(\h_0=0)=1\}
$$
\ie the infimum
of the values $q$ such that, with probability one, the lattice can be
entirely emptied.
Using the fact that the $c_x$'s are increasing function of $\h$ it is
easy to check that for any $q>q_{bp}$  $\mu^{(\infty)}(\h_0=0)=1$.
\begin{Proposition}
\label{qc=qbp} $q_c=q_{bp}$ and for any $q>q_c$
  $0$ is a simple eigenvalue for $\cL$.
\end{Proposition}
\begin{proof}
  Assume $q<q_{bp}$ and call $f$ the indicator of the event that the
  origin cannot be emptied by any finite number of iterations of the
  bootstrap map $T$ \eqref{eq:bootstrap map}. By construction
  $\Var(f)\neq 0$ and $\cL f=0$ a.s. $(\mu)$. Therefore $0$ is not a
  simple eigenvalue of $\cL$
  and $q\le q_c$.

  Suppose now that $q>q_{bp}$ and that $f\in Dom(\cL)$ satisfies $\cL
  f=0$ or, what is the same, $\cD(f)=0$. We want to conclude that
  $f=\text{const.}$ a.e. ($\mu$). For this purpose we will show that
  $\cD(f)=0$ implies that the unconstrained Glauber Dirichlet form
  $\sum_x\mu\left(\Var_x(f)\right)$ is zero which makes the sought
  conclusion obvious since $\Var(f)\le \sum_x
  \mu\left(\Var_x(f)\right)$.

  Given $x\in \bbZ^d$ let $A_n\equiv A_{n,x}=\{\h:\; T^n(\h)_x=0
  \}$. Since $q>q_{bp}$, clearly $\mu\left(\cup_n A_n\right)=1$. Write
  \begin{gather*}
\mu\left(\Var_x(f)\right)= pq\sum_n \int_{A_n\setminus A_{n-1}}d\mu(\h) [f(\h^x)-f(\h)]^2
  \end{gather*}
where $\h^x$ denotes the flipped configuration at $x$. For any $\h\in A_n$
it is easy to convince oneself that it is possible to find a collection
of vertices $x^{(1)},\dots,x^{(k)}$, with $k$ and $d(x,x^{(j)})$ bounded by a
constant depending only on $n$, and a collection of configurations
$\h^{(1)}, \h^{(2)},\dots,\h^{(k)}$ such that $\h^{(1)}=\h$,
$\h^{(k)}=\h^x$, $\h^{(j+1)}=(\h^{(j)})^{x^{(j)}}$ and
$c_{x^{(j)}}(\h^{(j)})=1$.
We can then write $[f(\h^x)-f(\h)]$ as a telescopic sum of terms like
$[f(\h^{(j+1)})-f(\h^{(j)})]$
and apply Schwartz inequality to get
\begin{gather*}
  \int_{A_n\setminus A_{n-1}}\!\! d\mu(\h) [f(\h^x)-f(\h)]^2 \\ \le
  C(n)\!\!\!\!\sum_{y:\,d(y,x)\le C'(n)}\int d\mu(\s)c_y(\s)[f(\s^y)-f(\s)]^2
\end{gather*}
where the constant $C(n)$ takes care of the relative density
$\sup_{\h\in A_n} \frac{\mu(\h)}{\mu(\h^{(j)})}$ and of the number
of possible choice of the
vertices $\{x^{(j)}\}_{j=1}^k$. \\
By assumption $\cD(f)=0$ \ie $\int
d\mu(\s)c_y(\s)[f(\s^y)-f(\s)]^2=0$ for any $y$ and the proof is complete.
\end{proof}
Having defined the bootstrap percolation it is natural to divide the 0-1
KCSM into two distinct classes.
\begin{definition}
We will say that a 0-1 KCSM is \emph{non cooperative} if there exists a
finite set $\cB\sset \bbZ^d$ such that any configuration $\h$ which is empty in all the
sites of $\cB$ reaches the empty configurations (all 0's) under
iteration of the bootstrap mapping. Otherwise the model will be called \emph{cooperative}.
\end{definition}
\begin{remark}
  Because of the translation invariance of the constraints it is obvious
  that any configuration $\h$ identically
equal to zero in $\cB +x$,
  $x\in\bbZ^d$, will reach the empty configuration under iterations of
  $T$. It is also obvious that $q_{bp}$ and therefore $q_c$ are zero for
  all non-cooperative models.
\end{remark}
In what follows we will now illustrate some of the most studied models.

\\
{\bf [1] The East model \cite{JE}.} Take $d=1$ and set $\cC_x=\{x+1\}$, \ie a
vertex can flip iff its right neighbor is empty. The minimal boundary
conditions in finite volume are of course empty right boundary. The
model is clearly cooperative but $q_c=0$ since in order to empty the
whole lattice it is enough to start from a configuration for which any
site $x$ has a vacancy to its right.

\\
{\bf [2] Frederickson-Andersen (FA-jf) models \cite{FA1,FA2}.} Take $1\le j\le d$ and set
$$
\cC_x=\{A\sset \cN_x: |A|\ge j\}
$$
In words a vertex can be updated iff at least $j$ of its neighbors are
$0$'s. When $j=1$ the minimal boundary conditions on a rectangle that will ensure
ergodicity of the Markov chain in e.g. a rectangle $\L$ will be
exactly one $0$ on $\partial \L$.  If instead $j=d$, ergodicity is
guaranteed if we assume $\t_y=0$ for e.g. all $y$ on $\partial_+\L$.  If
$j=1$ the model is non-cooperative while for $j\ge 2$ it is cooperative. In
any case $q_c=0$ \cite{Schonmann}.

\\
{\bf [3] The Modified Basic (MB) model.} Here we take
$$
\cC_x=\{A\sset \cN_x: A\cap\{-\vec e_i,\vec e_i\}\neq \emptyset, \text{
  for all $i=1,\dots,d$}\}
$$
\ie a move at $x$ can occur iff in each direction there is a $0$.
The model is cooperative and the minimal boundary conditions on a rectangle are the same as those for the
FA-df model. Once again $q_c=0$ \cite{Schonmann}.

\\
{\bf [4] The N-E (North-East) model \cite{RJM}.} Here one chooses $d=2$ and
$$
\cC_x=\{\cK_x\}
$$
The model is cooperative with minimal boundary conditions those that we
have chosen for the
FA-2f model in $d=2$. The critical point $q_c$ coincides with $1-p^o_c$
where $p_c^o$ is the critical threshold for oriented percolation in
$\bbZ^2$ \cite{Schonmann}.

\subsection{Quantities of interest}
Back to the  general model we now define the main quantities that will
be studied in the sequel.

The first object of mathematical and physical interest is the spectral
gap (or inverse of the relaxation time) of the generator $\cL$, defined as
\begin{equation}
  \label{eq:gap}
\gap(\cL):=\inf_{\substack{f\in Dom(\cL)\\ f\neq \text{const}}}\frac{\cD(f)}{\Var(f)}
\end{equation}
and similarly for the finite volume version of the process. A positive
spectral gap implies that the reversible measure $\mu$ is mixing for the
semigroup $P_t$ with exponentially decaying correlations.  It is
important to observe the following kind of monotonicity that can be
exploited in order to bound the spectral gap of one model with the
spectral gap of another one.

Suppose that we are given two finite range and translation invariant
influence classes $\cC'_0,\cC_0$ such that, for all $\o\in \O$ and all
$x\in \bbZ^d$,
$c_x(\o)\le c'_x(\o)$ and denote the associated generators by $\cL$ and $\cL'$
respectively. In this case we say that the KCSM generated by
${\mathcal{L}}$ is dominated by the one generated by ${\mathcal{L'}}$. Clearly $c_x(\o)\le c_x'(\o)$ for all $\o$ and
therefore $\gap(\cL)\le \gap(\cL')$.
As an example we can consider the FA-1f model in $\bbZ^d$. If instead of
taking as $\cC_0$ the collection of non-empty subsets $A$ of $\cN_0$
(see above) we consider  $\cC_0$ with the extra constraint that $A$ must
contain at least one vertex between $\pm \vec e_1$, we  get
that the spectral gap of the FA-1f model in $\bbZ^d$ is bounded from
below by the spectral gap of the FA-1f model in $\bbZ$ which in turn is
bounded from below by the spectral gap of the East model which is known
to be positive \cite{Aldous}.  Similarly we could lower bound
the spectral gap of the FA-2f model in $\bbZ^d$, $d\ge 2$, with that in
$\bbZ^2$, by restricting the sets $A\in \cC_0$ to e.g. the $(\vec e_1,
\vec e_2)$-plane.  In finite volume the comparison
argument is a bit more delicate since it heavily depends on the boundary
conditions. For example, if we consider the FA-1f model in a rectangle
with minimal boundary conditions, \ie a single $0$ in one
corner, the argument discussed above would lead to a comparison with a
non-ergodic Markov chain whose spectral gap is zero.
\begin{remark}
  The comparison technique can be
  quite effective in proving positivity of the spectral
  gap but the resulting bounds are in general quite poor,
  particularly in the limiting case $q\approx q_c$.
\end{remark}
The second observation we make consists in relating $\gap(\cL)$ to its
finite volume analogue.  Assume that $\inf_{\L\in\bbF}\gap(\cL_\L)>0$
where $\cL_\L$ is defined with e.g. good boundaries conditions outside
$\L$. It is then easy to conclude that $\gap(\cL)>0$.

Indeed, following Liggett Ch.4 \cite{Liggett}, for any $f\in Dom(\cL)$
with $\Var(f)>0$ pick $f_n\in L^2(\O,\mu)$ depending only on finitely
many spins so that $f_n \to f$ and $\cL f_n \to \cL f$ in $L^2$. Then
$\Var(f_n)\to \Var(f)$ and $\cD(f_n)\to \cD(f)$. But since $f_n$ depends
on finitely many spins
$$
\Var(f_n)=\Var_\L(f_n)\quad \text{and}\quad \cD(f_n)=\cD_\L(f_n)
$$
provided that $\L$
is a large enough square (depending on $f_n$) centered at the origin. Therefore
$$
\frac{\cD(f)}{\Var(f)}\ge \inf_{\L\in\bbF}\gap(\cL_\L)>0.
$$
and $\gap(\cL)\ge \inf_{\L\in\bbF}\gap(\cL_\L)>0$.

\medskip

The second quantity of interest is the so called \emph{persistence
  function} (see e.g. \cite{Ha,SE}) defined by
\begin{equation}
  \label{eq:Pers}
F(t):=\int d\mu(\h)\; \bbP(\s^\h_0(s)=\h_0,\; \forall s\le t)
\end{equation}
where $\{\s^\h_s\}_{s\ge 0}$ denotes the process started from the
configuration $\h$. In some sense the persistence function provides a
measure of the ``mobility'' of the system.

\section{Main Results for 0-1 KCSM}
\label{main results}
In this section we state our main results for 0-1 KCSM.
Fix an integer length scale $\ell$ larger than the range $r$ and let
$\bbZ^d(\ell)\equiv \ell\,\bbZ^d$. Consider a
partition of $\bbZ^d$ into disjoint rectangles $\L_z:=\L_0+z$, $z\in
\bbZ^d(\ell)$, where $\L_0=\{x\in\bbZ^d:\ 0\le x_i\le \ell-1,\; i=1,..,d\}$.
\begin{definition}
\label{defgelle}
Given $\epsilon\in (0,1)$ we say that $G_{\ell}\sset \{0,1\}^{\L_0}$ is
a $\epsilon$-good set of configurations on scale $\ell$ if the following
two conditions are satisfied:
\begin{enumerate}[(a)]
\item $\mu(G_{\ell})\ge 1-\epsilon$.
\item For any collection $\{\xi^{(x)}\}_{x\in \cK_0^*}$ of spin
  configurations in $\{0,1\}^{\L_0}$ such that $\xi^{(x)}\in G_\ell$ for
  all $x\in \cK_0^*$ and for any $\xi\in \O$
  which coincides with $\xi^{(x)}$ in $\L_{\ell x}$, there exists a
  sequence of legal moves inside $\cup_{x\in \cK_0^*}\L_{\ell x}$ (\ie
  single spin moves compatible with the constraints) which transforms
  $\xi$ into a new configuration $\t$ such that the Markov chain
  generated by $\cL_{\L_0}$ with boundary conditions $\t$ is ergodic.
\end{enumerate}
\end{definition}
\begin{remark} In general the transformed configuration $\tau$ will be
  identically equal to zero on $\partial_+^* \L_0$. It is also clear
  that assumption (b) has been made having in mind models,
  like the FA-jf, M-B or N-E, which, modulo rotations, are dominated by
  a model with influence classe  $\tilde \cC_x$ entirely
  contained in the sector $\{y: y=x+\sum_{i=1}^d \a_i \vec e_i,\ \a_i\ge 0\}$. Moreover,
  for some other models, the geometry of the tiles of the partition of
  $\bbZ^d$, rectangles in our case, should be adapted to the influence
  classes $\{\cC_x\}_{x\in \bbZ^d}$.
\end{remark}
% \begin{remarks} \
%   \begin{enumerate}[i)]
% \item In general the transformed configuration $\tau$ will be
%   identically equal to zero on $\partial_+^* \L_0$
% \item For some specific models, like e.g. the Knights
%   models introduced in \cite{Toninelli}, the geometry of the tiles of the partition
%   of $\bbZ^d$, rectangles in our case, should be adapted to the
%   influence class $\{\cC_x\}_{x\in \bbZ^d}$.
%   \end{enumerate}
% \end{remarks}
With the above notation our first main result, whose proof can be found
in section \ref{Proof of theorem 1}, can be formulated as follows.
\begin{Theorem}
There exists a universal constant $\epsilon_0\in (0,1)$ such that if
there exists $\ell$ and a $\epsilon_0$-good set $G_\ell$ on scale $\ell$
then $\gap(\cL)>0$.
\label{main theorem 01}
\end{Theorem}
In several examples, e.g. the FA-jf and Modified
  Basic models, the natural candidate for the event
  $G_\ell$ is  the event that the tile $\L_0$ is ``internally
    spanned'', a notion borrowed from bootstrap percolation
  \cite{Aizenman,Schonmann,Holroyd}:
  \begin{definition}
\label{ISDEF}
We say that a finite set $\G\subset\bbZ^d$ is \emph{internally spanned}
by a configuration $\eta\in\O$ if, starting from
the configuration $\h^\G$ equal to one outside $\G$ and equal to $\h$
inside $\G$,  there exists a sequence of legal moves
  inside $\G$ which connects
  $\h^\G$ to the configuration identically equal to zero inside $\G$ and
  identically equal to one outside $\G$.
  \end{definition}
 Of course whether or not the
  set $\L_0$ is internally spanned for $\h$ depends only on the
  restriction of $\h$ to $\L_0$. One of the major result in
  bootstrap percolation problems has been the exact evaluation of the
  $\mu$-probability that the box $\L_0$ is internally spanned as a
  function of the length scale $\ell$ and the parameter $q$
  \cite{Holroyd, Schonmann, Cerf, Aizenman}. For non-cooperative models
  it is obvious that $q_{bp}=0$.
  For some cooperative systems like
  e.g. the FA-2f and Modified Basic model in $\bbZ^2$, it has been shown that for any
  $q>0$ such
  probability tends very rapidly (exponentially fast) to one as $\ell\to
  \infty$ and that it abruptly jumps from being very small to being
  close to one as $\ell$ crosses a critical scale
  $\ell_c(q)$. In most cases the critical length $\ell_c(q)$
  diverges very rapidly as $q\downarrow 0$. Therefore, for such models
  and $\ell>\ell_c(q)$, one could safely take $G_\ell$ as the collection of
  configurations $\h$ such that
  $\L_0$ is internally spanned for $\h$. We now formalize what we just
  said.
  \begin{Corollary}
\label{IS}
  Assume that $\lim_{\ell\to \infty}\mu(\L_0 \text{
    is internally spanned })=1$ and that the Markov chain in $\L_0$ with
  zero boundary conditions on $\cup_{x\in \cK_0^*}\L_{\ell x}$ is ergodic. Then $\gap(\cL)>0$.
  \end{Corollary}
The second main result concerns the long time behavior of the
persistence function $F(t)$ defined in \eqref{eq:Pers}.

\begin{Theorem}
\label{persf}
Assume that $\gap(\cL)\ge \g>0$. Then there exists a constant $c>0$ such
that $F(t)\le \nep{-ct}$. For small values of $\g$ the constant $c$ can
be taken
proportional to $q\g$.
\end{Theorem}
\begin{proof}
Clearly
$F(t)=F_1(t)+F_0(t)$ where
$$
F_1(t)= \int\, d\mu(\h)\,
\bbP(\s^\h_0(s)=1\ \text{for all $s\le t$})
$$
and similarly for $F_0(t)$.  We will prove the exponential decay of
$F_1(t)$ the case of $F_0(t)$ being similar.

For any $\l>0$ the exponential Chebychev inequality gives
$$
F_1(t)= \int\, d\mu(\h)\,
\bbP\Bigl(\int_0^t ds\, \s^\h_0(s)=t\Bigr)
\le \nep{-\l t}\ \bbE_\mu\bigl(\nep{\l\int_0^t ds \,\s^\h_0(s)}\bigr)
$$
where $\bbE_\mu$ denotes the expectation over the process started from
the equilibrium distribution $\mu$. On $L^2(\mu)$
consider the self-adjoint operator $H_\l:=\cL + \l V$, where
$V$ is the multiplication operator by $\s_0$.
By the very definition of the scalar product $\scalprod{f}{g}$ in
$L^2(\mu)$ and the Feynman--Kac formula, we can rewrite $\bbE_\mu(\nep{\l\int_0^t \s_0(s)})$ as
$\scalprod{1}{\nep{tH_\l}1}$.
Thus, if $\b_\l$ denotes the supremum of the spectrum of $H_\l$,
$$
\bbE_\mu(\nep{\l\int_0^t \s_0(s)})\leq \nep{t\b_\l}.
$$
In order to complete the proof we need to show that for suitable
positive $\l$ the constant $\b_\l/\l$ is strictly smaller than one.

For any norm one function $f$ in the domain of $H_\l$ (which coincides
with Dom($\cL$)) write $f=\a {\bf 1} +g$ with
$\scalprod{1}{g}=0$. Thus
\begin{gather}
\scalprod{f}{H_\l f}=\scalprod{g}{\cL g} \! + \a^2 \l\scalprod{1}{V
1} \!\! +
\l\scalprod{g}{Vg} \! \!+ 2\l\a\scalprod{1}{Vg} \nonumber \\
\le (\l-\g) \scalprod{g}{g} +\a^2\l p +
2\l|\a|\bigl(\scalprod{g}{g}pq\bigr)^{1/2}
\end{gather}
Since $\a^2+\scalprod{g}{g}=1$
\begin{equation}
  \b_\l/\l \le \sup_{0\le \a\le 1}\Bigl\{(1-\g/\l)(1-\a^2)+p\a^2+2\a\bigl((1-\a^2)pq\bigr)^{1/2}\Bigr\}
\label{eq:pers2}
\end{equation}
If we choose $\l=\g/2$ the r.h.s. of \eqref{eq:pers2}
becomes
\begin{gather*}
\sup_{0\le \a\le 1} (1+p)\a^2 -1 + 2\a\bigl((1-\a^2)pq\bigr)^{1/2} \\
\le \sup_{0\le \a\le 1} (1+p)\a^2 -1 + 2\bigl((1-\a^2)pq\bigr)^{1/2}=
\frac{pq}{1+p}+p <1 .
\end{gather*}
since $p\neq 1$.
Thus $F_1(t)$ satisfies
$$
F_1(t)\le \nep{-t\frac{\g}{2} \frac{q}{1+p}}.
$$
A similar computation shows that $F_0(t)\le \nep{-t\g c}$
with $c$ independent of $q$.
\end{proof}
\begin{remark}
The above result indicates that one can obtain \emph{upper bounds} on
the spectral gap by proving \emph{lower bounds} on the persistence
function. Concretely a lower bound on the persistence function can be obtained
by restricting the $\mu$-average to those initial configurations
$\h$ for which the origin is blocked with high probability for all times
$s\le t$.  In section \ref{specificmodels}  we will see few examples of
this strategy.
\end{remark}

\section{Analysis of a general auxiliary model}
\label{GM}
Consider the following model characterized by the influence classes
$\cC_x=\cK_x^*,\ x\in \bbZ^d$ and arbitrary
finite probability space $(S,\nu)$ and choice of the good event
$G\sset S$ with $q:=\nu(G)$. For definiteness we will call it the
\emph{*-general model}. The proof of theorem \ref{main theorem 01} is based on the analysis of
the *-general model in a finite set $\L$
with fixed \emph{good} boundary conditions $\t$ on its *-oriented
neighborhood $\partial_+^*\L$. Clearly the process does not depend on the
specific values of the (good) boundary configuration $\t$ and, with a
slightly abuse of notation, we can safely denote
the generator of the chain by $\cL_\L$ and the associated Dirichlet form
by $\cD_\L$. Ergodicity of $\cL_\L$ follows  once we
observe that, starting from the sites in $\L$ whose *-oriented
neighborhood is entirely contained in $\L^c$ and whose existence is
 proved by induction, we can reach any good
configuration $\o'\in G^\L$ and from there any other configuration
$\tilde\o$.

The following monotonicity of the spectral gap will turn out to be quite useful in simplifying
some of the arguments given below.
\begin{Lemma}
\label{Mon}
Let $V\sset \L$. Then
$$
\gap(\cL_\L)\le \gap(\cL_V)
$$
\end{Lemma}
\begin{proof}
For any $f\in L^2(\O_V,\mu_V)$ we have $\Var_V(f)=\Var_\L(f)$ because
of the product structure of the measure $\mu_\L$ and
$\cD_\L(f)\le \cD_V(f)$ because, for any $x\in V$ and any $\o\in
\O_{\L}$, $c_{x,\L}(\o)\le c_{x,V}(\o)$. The
result follows at once from the variational characterization of the
spectral gap.
\end{proof}
We now state our main theorem concerning the *-general model.
\begin{Theorem}
\label{main theorem}
There exists $q_0<1$ independent of $S,\nu$ such that for any $q>q_0$
$$
\inf_{\L\in \bbF } \gap(\cL_\L)>1/2.
$$
and in particular $\gap(\cL)>0$.
\end{Theorem}
\begin{proof}
Thanks to Lemma \eqref{Mon} we need to prove the result only for
rectangles.
Our approach is based on the ``bisection method'' introduced in
\cite{SFlour,Encyclopedia} and which, in its essence, consists in
proving a suitable recursion relation between the spectral gap on scale
$2L$ with that on scale $L$. At the beginning the method requires a
simple geometric result (see \cite{Cesi}) which we now describe.

Let $l_k := (3/2)^{k/2}$, and let $\bbF_k$ be the set of all rectangles
$\L\sset \bbZ^d$
which, modulo translations and permutations of the coordinates, are
contained in
$$
  [ 0,l_{k+1} ] \times \dots\times [0,l_{k+d}]
$$
The main property of $\bbF_k$ is that each rectangle in
$\bbF_k\setm \bbF_{k-1}$ can be obtained as a ``slightly
overlapping union'' of two rectangles in $\bbF_{k-1}$. More precisely
we have:

\begin{Lemma}
\label{geom}
For all $k\in \bbZ_+$,
for all $\L \in \bbF_k\setm \bbF_{k-1}$ there exists a finite sequence
$\{\L_1^{(i)}, \L_2^{(i)}\}_{i=1}^{s_k}$ in $\bbF_{k-1}$, where
$s_k := \inte{l_k^{1/3}}$, such that, letting $\d_k := \frac 18
\sqrt{l_k}-2$,
\begin{enumerate}[(i)]
\item $\L = \L_1^{(i)} \cup \L_2^{(i)}$,
\item  $d(\L\setm \L_1^{(i)}, \L\setm \L_2^{(i)}) \ge \d_k $,
\item $\left(\L_1^{(i)}\cap \L_2^{(i)}\right)\cap \left(\L_1^{(j)}\cap \L_2^{(j)}\right)
        = \emp$,  if $i\ne j$
\end{enumerate}
\end{Lemma}
The bisection method then establishes a simple recursive inequality
between the quantity $\g_k := \sup_{\L\in\bbF_k} \gap(\cL_\L)^{-1}$ on
scale $k$ and the same quantity on scale $k-1$ as follows.

Fix $\L\in \bbF_k\setm \bbF_{k-1}$ and write it as $\L = \L_1 \cup \L_2$ with $\L_1,
\L_2\in \bbF_{k-1}$ satisfying the properties described in Lemma \ref{geom}
above. Without loss of generality we can assume that
all the faces of $\L_1$ and of $\L_2$ lay on the faces of
$\L$ except for one face orthogonal to the first direction $\vec e_1$
and that, along that direction, $\L_1$ comes before $\L_2$.
%(see figure \eqref{}).
Set $I\equiv \L_1 \cap \L_2$ and write, for definiteness,
$I=[a_1,b_1]\times\dots\times[a_d,b_d]$. Lemma \ref{geom} implies that
the width of $I$ in the first direction, $b_1-a_1$, is at least
$\d_k$.
Let also
% $I_r=[(a_1+b_1)/2,b_1]\times\dots\times[a_d,b_d]$,
% $I_l=I\setminus I_r$ and
$\partial_{r}I=\{b_1\}\times\dots\times[a_d,b_d]$
% be the right/left half and
be the right face of $I$ along the first direction.

Next, for any $x,y\in I$ and any $\o\in \O_{I}$, we write $x
\stackrel{\o}{\rightarrow}y$ if there exists a sequence
$(x^{(1)},\dots,x^{(n)})$ in $I$, starting at $x$ and ending at $y$,
such that, for any $j=1,\dots, n-1$, $x^{(j)}\sim x^{(j+1)}$ and
$\o_{x^{(j)}}\notin G$, where $\sim$ has been defined in section
\ref{setting}.  With this notation we finally define the \emph{bad
  cluster} of $x$ as the set $A_x(\o)=\{ y\in I;\
x\stackrel{\o}{\rightarrow} y\}$. Notice that, by construction, $\o_z\in
G$ for any $z\in \partial^* A_x(\o)$.

\begin{definition}
 We will say that $\o$ is $I$-good
  iff, for all $x\in \partial_{r}I$, the set $A_x(\o)\cup\partial^*A_x(\o)$
  is contained in $I$.
\end{definition}

With the help of the above decomposition
we now run the following constrained ``block dynamics'' on $\O_\L$ (in what follows,
for simplicity, we suppress the index $i$) with blocks
$B_1:=\L\setminus \L_2$ and $B_2:=\L_2$. The block $B_2$ waits a mean one
exponential random time and then the current configuration inside it is refreshed
with a new one sampled from $\mu_{\L_2}$. The block $B_1$ does the same
but now the configuration is refreshed only if the current configuration
$\o$ is $I$-good (see Figure \ref{Fig:blocchi}).
\begin{figure}[h]
\psfrag{1}{$B_1$}
\psfrag{2}{$B_2$}
\psfrag{3}{$I$}
\psfrag{6}{$\partial_r I$}
\includegraphics[width=2.7in]{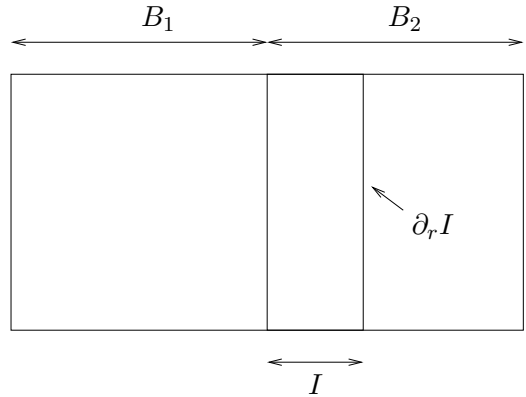}
\caption{The two blocks and the strip $I$.}
\label{Fig:blocchi}
\end{figure}
The Dirichlet form of this auxiliary
chain is simply
$$
\cD_{block}(f)=\mu_\L\left(c_1\Var_{B_1}(f)+\Var_{B_2}(f)\right)
$$
where $c_1(\o)$ is just the indicator of the event
that $\o$ is $I$-good and
$\Var_{B_1}(f)$, $\Var_{B_2}(f)$ depend on $\o_{B_1^c}$ and
$\o_{B_2^c}$ respectively.

Denote by $\g_{\rm block}(\L)$ the inverse spectral gap of this auxiliary
chain. The following bound, whose proof is postponed for clarity of the
exposition, is not difficult to prove.

\begin{Proposition}
\label{gapblock}
Let $\e_k\equiv {\displaystyle \max_{I}}\,\bbP(\o\text{ is not
  $I$-good})$ where the $\max_{I}$ is taken over the $s_k$ possible choices of
the pair $\left(\L_1,\L_2\right)$. Then
$$
\g_{\rm block}(\L) \le \frac{1}{1-\sqrt{\e_k}}
$$
\end{Proposition}
In conclusion, by writing down the standard Poincar\'e inequality for the
block auxiliary chain, we get that for any $f$
\begin{equation}
  \label{eq:s1}
\Var_\L(f)\le \bigl(\frac{1}{1-\sqrt{\e_k}}\bigr)\ \mu_\L\Bigl(c_1\Var_{B_1}(f)+\Var_{B_2}(f)\Bigr)
\end{equation}

The second term, using the definition of $\g_k$ and the fact that
$B_2\in \bbF_{k-1}$ is bounded from above
by
\begin{equation}
  \label{eq:s2}
  \mu_\L\Bigl(\Var_{B_2}(f)\Bigr)\le \g_{k-1}\sum_{x\in
    B_2}\mu_\L\bigl(c_{x,B_2} \Var_x(f)\bigr)
\end{equation}
Notice that, by construction, for all $x\in B_2$ and all $\o$,
$c_{x,B_2}(\o)=c_{x,\L}(\o)$.  Therefore the term $\sum_{x\in
  B_2}\mu_\L\bigl(c_{x,B_2} \Var_x(f)\bigr)$ is nothing but the
contribution carried by the set $B_2$ to the full Dirichlet form $\cD_\L(f)$.

Next we examine the more complicate term
$\mu_\L\Bigl(c_1\Var_{B_1}(f)\Bigr)$ with the goal
in mind to bound it with the missing term of the  full Dirichlet form
$\cD_\L(f)$.

For any $I$-good $\o$ let
$\Pi_\o=\cup_{x\in \partial_{r}I}A_x(\o)$, let $B_\o$ be the
connected (w.r.t. the graph structure induced by the $\sim$
relationship)
component of $B_1\cup I\setminus \left(\Pi_\o \cup \partial^*
  \Pi_\o\cup\partial_{r}I\right)$ which contains $B_1$ (see Figure \ref{Fig:piomegabis}).
\begin{figure}[h]
\psfrag{1}{$B_1$}
\psfrag{2}{$I$}
\psfrag{3}{$\partial_+ \Lambda$}
\psfrag{4}{$\partial_r I$}
\includegraphics[height=3in]{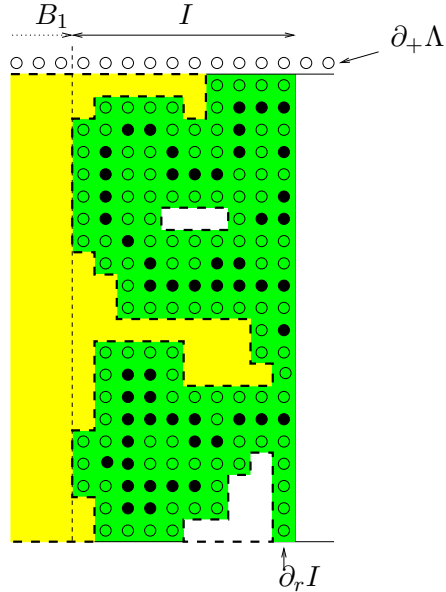}
\caption{An example of an $I$-good configuration $\omega$: empty sites are
good and filled ones are noT good. The grey region is the set $\Pi_\omega \cup
\partial^* \Pi_\omega \cup \partial I_r$. The dotted lines mark the
connected components of
$B_1 \cup I \setminus (\Pi_\omega \cup \partial^* \Pi_\omega \cup \partial
I_r)$. The connected component containing $B_1$ is the shaded one.}
\label{Fig:piomegabis}
\end{figure}
A first key observation is
now the following.
\begin{claim}
\label{claim1}
For any $z\in \partial_+^*B_\o$ it holds true that  $\o_z\in G$ .
\end{claim}
\begin{proof}[Proof of the claim]
  To prove the claim suppose the opposite and let
  $z\in \partial_+^*B_\o$ be such that $\o_z\notin G$ and let $x\in
  B_\o$ be such that $\cK_x^*\ni z$. Necessarily $z\in \Pi_\o$ because
  of the good boundary conditions in $\partial_+^*\L$ and the fact that
  $\o_y\in G$ for all $y\in \partial^*\Pi_\o\cup
  \left(\partial_rI\setminus \Pi_\o\right)$. However $z\in \Pi_\o$ is
  impossible because in that case $z\in A_y(\o)$ for some
  $y\in \partial_{r}I$ and therefore $x\in A_y(\o)\cup\partial^*
  A_y(\o)$ \ie $x\in \Pi_\o\cup\partial^*\Pi_\o$, a contradiction.
\end{proof}
The second observation is the following.
\begin{claim}
\label{claim2}
  For any $\G\sset \O_\Pi:=\cup_{\o\ \text{I-good}}\Pi_\o$, the
  event $\{\o:\ \Pi_\o=\G\}$ does not depends on the values of $\o$ in
  $B_\G$, the connected component (w.r.t. $\sim$) of\\ $B_1\cup I\setminus
  \G\cup\partial^*\G\cup \partial_rI$ which contains $B_1$.
\end{claim}
\begin{proof}[Proof of the claim]
Fix $\G\in \O_\Pi$. The event $\Pi_\o=\G$ is equivalent to:
\begin{enumerate}[(i)]
\item $\o_z\in G$ for any
  $z\in \partial_rI\setminus\G$;
\item $\o_z\in G$ for any
  $z\in \partial^*\G\cap I$;
\item $\o_z\notin G$ for all $z\in \G$.
\end{enumerate}
In fact trivially $\Pi_\o=\G$ implies (i),(ii) and (iii). To prove the
other direction we first observe that (i) and (iii) imply that
$\Pi_\o\supset \G$. If $\Pi_\o\neq \G$ there exists $
z\in\Pi_\o\setminus\G$ which is in $\partial^*\G\cap I$ and such that
$\o_z\notin G$. That is clearly impossible because of (ii).
\end{proof}
If we
observe that $\Var_{B_1}(f)$ depends only on $\o_{B_2}$, we can write
(we omit the subscript $\L$ for simplicity)
\begin{gather}
  \mu\Bigl(c_1\Var_{B_1}(f)\Bigr)=\sum_{\G\in \O_\Pi}\mu\bigl(\id_{\{\Pi_\o=\G\}}\Var_{B_1}(f)\bigr)
  \nonumber \\
=\sum_{\G\in \O_\Pi}\sum_{\o_{B_2\setminus I}}\mu(\o_{B_2\setminus
  I})\sum_{\o_I}\mu(\o_I)\id_{\{\Pi_\o=\G\}}\Var_{B_1}(f)\nonumber \\
=\sum_{\G\in \O_\Pi}\sum_{\o_{B_2\setminus I}}\mu(\o_{B_2\setminus I})
\sum_{\o_{I\setminus I_\G}}\mu(\o_{I\setminus I_\G})\id_{\{\Pi_\o=\G\}}\sum_{\o_{I_\G}}\mu(\o_{I_\G})\Var_{B_1}(f)
\label{A}
\end{gather}
where $I_\G=B_\G \cap I$ and we used the independence
of $\id_{\{\Pi_\o=\G\}}$ from $\o_{I_\G}$.

The convexity of the variance implies that
\begin{equation*}
\sum_{\o_{I_\G}}  \mu(\o_{I_\G})\Var_{B_1}(f)\le \Var_{B_\G}(f)
\end{equation*}
The Poincar\'e inequality together
with Lemma \eqref{Mon} finally gives
\begin{gather}
\Var_{B_\G}(f)\le \gap(\cL_{B_\G})^{-1}\sum_{x\in
  B_\G}\mu_{B_\G}\bigl(c_{x,B_\G}\Var_{x}(f)\bigr)\nonumber \\
\le \gap(\cL_{B_1\cup I})^{-1}\sum_{x\in
  B_\G}\mu_{B_\G}\bigl(c_{x,B_\G}\Var_{x}(f)\bigr)
\label{B}
\end{gather}
The role of the event
$\{\Pi_\o=\G\}$ should at this point be clear. For any $\o\in\O_\L$ such that
$\Pi_\o=\G$, let $\o_{B_\G}$ be its restriction to the set
$B_\G$. From claim \ref{claim1} we  infer that
\begin{equation}
  \label{C}
c_{x,\L}(\o)=c_{x,B_\G}(\o_{B_\G})\quad \forall x\in B_\G\,.
\end{equation}
If we finally plug \eqref{B} and \eqref{C} in the r.h.s. of \eqref{A}
and recall that $B_1\cup I=\L_1\in \cF_{k-1}$,
we obtain
\begin{gather}
  \mu_\L\Bigl(c_1\Var_{B_1}(f)\Bigr)\le
\gap(\cL_{\L_1})^{-1}\mu_\L\bigl(c_1\sum_{x\in
  B_{\Pi_\o}}c_{x,\L}\Var_{x}(f)\bigr)\nonumber \\
\le
\g_{k-1}\,\mu_\L\bigl(\sum_{x\in
  \L_1}c_{x,\L}\Var_{x}(f)\bigr)
\label{D}
\end{gather}

In conclusion we have shown that
\begin{equation}
  \label{eq:2}
  \Var_\L(f)\le \bigl(\frac{1}{1-\sqrt{\e_k}}\bigr)\g_{k-1}\Bigl(\cD_\L(f)+\sum_{x\in \L_1\cap\L_2}\mu_\L\bigl(c_{x,\L}\Var_x(f)\bigr)\Bigr)
\end{equation}
Averaging over the $s_k= \inte{l_k^{1/3}}$ possible choices of
the sets $\L_1,\L_2$ gives
\begin{equation}
  \label{eq:3}
  \Var_\L(f)\le \bigl(\frac{1}{1-\sqrt{\e_k}}\bigr)\g_{k-1}(1+\frac{1}{s_k})\cD_\L(f)
\end{equation}
which implies that
\begin{gather}
  \label{eq:4}
  \g_k\le \bigl(\frac{1}{1-\sqrt{\e_k}}\bigr)(1+\frac{1}{s_k})\g_{k-1}\\
\le \g_{k_0}\ \prod_{j=k_0}^k \bigl(\frac{1}{1-\sqrt{\e_j}}\bigr)(1+\frac{1}{s_j})
\end{gather}
where $k_0$ is the smallest integer such that $\d_{k_0}>1$.

It is at this stage (and only here) that we need a restriction on
the probability $q$ of the good set $G$. If $q$ is taken large enough (but uniformly
in the cardinality of $S$), the quantity $\e_j$ becomes exponentially
small in $\d_j=\frac 18 \sqrt{l_j}-2$ (the minimum width of the
intersection between the rectangles $\L_1,\L_2$ on scale $l_j$) with a
large constant rate and the convergence of the infinite product $
\prod_{j=k_0}^\infty
\bigl(\frac{1}{1-\sqrt{\e_j}}\bigr)(1+\frac{1}{s_j}) $ as well as the
fact that the quantity $\g_{k_0}\ \prod_{j=k_0}^k
\bigl(\frac{1}{1-\sqrt{\e_j}}\bigr)(1+\frac{1}{s_j})$ is smaller than
$2$ follows at once from the exponential growth of the scales
$l_j=(3/2)^{j/2}$.
\end{proof}

\begin{proof}[Proof of Proposition \eqref{gapblock}]
For any mean zero function $f\in L^2(\O_\L,\mu_\L)$ let
$$
\pi_1f:=\mu_{B_2}(f), \quad \pi_2f:=\mu_{B_1}(f)
$$
be the natural projections onto $L^2(\O_{B_i},\mu_{B_i})$,
$i=1,2$. Obviously $\pi_1\pi_2f=\pi_2\pi_1f=0$. The generator of the block dynamics can then be
written as:
$$
\cL_{\rm block}f=c_1\bigl(\pi_2f -f \bigr) +\pi_1f -f
$$
and the associated eigenvalue equation as
\begin{equation}
  \label{eq:4bis}
c_1\bigl(\pi_2f -f \bigr) +\pi_1f -f=\l f.
\end{equation}
By taking $f(\s_\L)=g(\s_{B_2})$ we see that $\l=-1$ is an eigenvalue. Moreover, since $c_1\le
1$, $\l\ge -1$.
Assume now $0>\l> -1$ and
apply $\pi_2$ to both sides of \eqref{eq:4bis} to obtain (recall that $c_1=c_1(\s_{B_2})$)
\begin{equation}
  \label{eq:5}
  -\pi_2f=\l\pi_2f \quad \imp \quad \pi_2f=0
\end{equation}
For any $f$ with $\pi_2f=0$ the eigenvalue equation becomes
\begin{equation}
  \label{eq:6}
  f=\frac{\pi_1 f}{1+\l+c_1}
\end{equation}
and that is possible only if
$$
\mu_{B_2}(\frac{1}{1+\l+c_1})=1.
$$
We can  solve the equation to get
$$
\l=-1+\sqrt{1-\mu_{B_2}(c_1)}\le -1+\sqrt{\e_k}.
$$
\end{proof}

\section{Proof of Theorem \ref{main theorem 01}}
\label{Proof of theorem 1}
In this section we provide the proof of the main Theorem \ref{main
  theorem 01}. For the relevant notation we refer the reader to section
\ref{main results}.

 Define $\epsilon_0=1-q_0$ where $q_0$ is the threshold appearing in
  Theorem \ref{main theorem} and assume that $\ell$ is such that there exists a $\epsilon_0$-good event
  $G_\ell$ on scale $\ell$.
  Consider the *-general model on $\bbZ^d(\ell)$ with $S=\{0,1\}^{\L_0}$,
  $\nu=\mu_{\L_0}$ and good event $G_\ell$. Obviously the two
  probability spaces $\O=\bigl(\{0,1\}^{\bbZ^d},\mu\bigr)$ and $
  \O(\ell)=\bigl(S^{\bbZ^d(\ell)},\prod_{x\in\bbZ^d(\ell)}\nu_x \bigr)$ coincide.
  Thanks to condition (a) on $G_\ell$ we can use theorem \ref{main
    theorem} to get that for any $f\in Dom(\cL)$
  \begin{equation}
    \Var(f)\le 2\sum_{x\in \bbZ^d(\ell)}\mu\bigl(\tilde c_x\Var_{\L_x}(f)\bigr)
 \label{fundam}
 \end{equation}
where the (renormalized) rate $\tilde c_x(\s)$ is simply the indicator
function of the event
that for any $y\in \cK^*_{\{x/\ell\}}$ the restriction of $\s$ to the
rectangle $\L_{\ell y}$ belongs to the good set $G_\ell$ on scale $\ell$.

In the sequel we will often refer to \eqref{fundam} as the
\emph{renormalized-Poincar\'e inequality} with parameters $(\ell,G_\ell)$.

Let us examine a generic term $  \mu\bigl(\tilde
c_x(\xi)\Var_{\L_x}(f)\bigr)$ which we write as
\begin{gather}
  \label{mt-1}
\frac 12 \int d\mu(\xi)
  \tilde c_x(\xi)\int \int d\mu_{\L_x}(\s)d\mu_{\L_x}(\h)
  \bigl[f(\s\cdot\xi)-f(\h\cdot\xi)\bigr]^2
\end{gather}
By assumption, if $\tilde c_x(\xi)=1$ necessarily there exists $\t$ and
a sequence of configurations $(\xi^{(0)},\xi^{(1)},\dots,\xi^{(n)})$, $n\le 3\ell^d$,
with the following properties:
\begin{enumerate}[(i)]
\item $\xi^{(0)}=\xi$ and $\xi^{(n)}=\t$;
\item the chain in $\L_x$ with boundary conditions $\tau$ is ergodic;
\item $\xi^{(i+1)}$ is obtained from $\xi^{(i)}$ by changing exactly only one
  spin at a suitable site $x^{(i)}\in \cup_{y\in \cK^*_{\{x/\ell\}}}\L_{\ell y}$;
\item the move at $x^{(i)}$ leading from $\xi^{(i)}$ to $\xi^{(i+1)}$ is
  permitted \ie $c_{x^{(i)}}(\xi^{(i)})=1$ for every $i=0,\dots,n$.
\end{enumerate}
\begin{remark}
Notice that for any $i=0,\dots,n$, the intermediate
configuration $\xi^{(i)}$ coincides with $\xi$ outside $\cup_{y\in
  \cK^*_{\{x/\ell\}}}\L_{\ell y}$. Therefore, given
  $\xi^{(i)}=\h$, the number of starting configurations $\xi=\xi^{(0)}$
  compatible with $\h$ is bounded from above by $2^{3\ell^d}$ and the
  relative probability $\mu(\xi)/\mu(\h)$ by $\bigl(\min(p,q)\bigr)^{3\ell^d}$.
\end{remark}
By adding and subtracting the terms $f(\s\cdot \tau), f(\h\cdot\tau)$
inside $ \bigl[f(\s\cdot\xi)-f(\h\cdot\xi)\bigr]^2$ and by writing $
f(\s\cdot\tau)-f(\s\cdot\xi)$ as a telescopic sum $\sum_{i=1}^{n-1}
\bigl[f(\s\cdot\xi_{i+1})-f(\s\cdot\xi_i)\bigr]$ we get
\begin{gather}
 \label{mt-2}
\bigl[f(\s\cdot\xi)-f(\h\cdot\xi)\bigr]^2
\le
3\bigl[f(\s\cdot\tau)-f(\h\cdot\tau)\bigr]^2
 \nonumber \\
\label{pippo}
+3n\sum_{i=1}^{n-1}
\bigl[f(\s\cdot\xi^{(i+1)})-f(\s\cdot\xi^{(i)})\bigr]^2
+3n\sum_{i=1}^{n-1} \bigl[f(\h\cdot\xi^{(i+1)})-f(\h\cdot\xi^{(i)})\bigr]^2
\end{gather}
If we plug \eqref{pippo} inside the r.h.s. of \eqref{mt-1} and use
properties (i),...,(iv) of the intermediate configurations
$\{\xi^{(i)}\}_{i=1}^n$ together with the remark and the fact that the inverse spectral gap
in $\L_x$ with ergodic boundary conditions $\tau$ is bounded from above
by a constant depending only on $(q,\ell)$, we  get that
there exists a finite constant $c:=c(q,\ell)$ such that
$$
\mu\bigl(\tilde c_x(\xi)\Var_{\L_x}(f)\bigr)\le
c\!\!\!\sum_{y\in \L_x\cup_{y\in \cK^*_{\{x/\ell\}}\L_{\ell y}}}\mu\bigl(c_y
\Var_y(f)\bigr)
$$
and the proof is complete.

\section{Specific models}
\label{specificmodels}
In this section we analyze the specific models that have been introduced
in section 2 and for each of them we prove positivity of the spectral
gap for $q>q_c$ together with upper and lower bounds bounds as
$q\downarrow q_c$.
\subsection{The East model}
As a first application of our bisection method we reprove the result
contained in
\cite{Aldous} on the positivity of the spectral gap, but we sharpen (by a power
of $2$) their
lower bound.
\begin{Theorem}
For any $q\in(0,1)$ the spectral gap of the East model is
positive. Moreover, for any $\d\in (0,1)$ there exists $C_\d>0$ such that
\begin{equation}
  \label{eq:th1}
\gap\ge C_\d q^{\log_2(1/q)/(2-\d)}
\end{equation}
In particular
\begin{equation}
\lim_{q\to 0} \log(1/\gap)/(\log(1/q))^2 =\left(2\log 2\right)^{-1}
\label{eq:th2}
\end{equation}
\end{Theorem}
\begin{remark}
Notice that \eqref{eq:th2} disproves the asymptotic behavior of
the spectral gap suggested in \cite{SE}.
\end{remark}
\begin{proof}
The limiting result \eqref{eq:th2} follows at once from the lower bound together with the
analogous upper bound proved in \cite{Aldous}.

In order to get the lower bound \eqref{eq:th1} we want to apply
directly the bisection method used in the proof of theorem
\ref{main theorem} but we need to choose the length scales $l_k$ a
little bit
more carefully.

Fix $\d\in (0,1)$ and define $l_k=2^k,\,
\d_k=\inte{l_k^{1-\d/2}}, s_k:=\inte{l_k^{\d/6}}$. Let also $\bbF_k$ be the set of
intervals which, modulo translations, have the form $[0,\ell]$ with
$\ell\in [l_k,l_k+l_k^{1-\d/6}]$ and define $\g_k$ as the worst case over the
elements $\L\in\bbF_k$ of the inverse spectral gap in
$\L$ with empty boundary condition at the right boundary of $\L$. Thanks
to lemma \ref{Mon} the worst case is attained for the interval
$\L_k=[0,l_k+l_k^{1-\d/6}]$. With these notation there exists $k_\d$
independent of $q$ such that the same result of
lemma \ref{geom} holds true as long as $k\ge k_\d$. We can then repeat exactly the same analysis
done in the proof of theorem \ref{main theorem} to get that
\begin{equation}
  \g_k\le \g_{k_\d}\prod_{j=k_\d}^\infty\left(\frac{1}{1-\sqrt{\e_j}}\right)\,\prod_{j=k_\d}^\infty\left(1+\frac{1}{s_j}\right)
\label{eq:east4}
\end{equation}
Here the
quantity $\e_k$  is just the
probability that an interval of width $\d_k$ is fully occupied
(see proposition \ref{gapblock}) \ie
$\e_k=p^{\d_k}$. The convergence of the product in (\ref{eq:east4}) is thus
guaranteed and the positivity of the spectral gap follows.

Let us now discuss the asymptotic behavior of the gap as $q\downarrow
0$. We first observe that $\g_{k_\d}< (1/q)^{\a_\d}$ for some finite
$\a_\d$. That follows e.g. from a coupling argument. In a time lag one
and  with
probability larger than $q^{\a_\d}$ for suitable $\a_d$,
any configuration in $\L_{k_\d}$ can reach the empty configuration by just flipping one
after another the spins starting from the right boundary. In other
words, under the maximal coupling,
two arbitrary configurations will couple in a
time lag one with probability larger than $q^{\a_\d}$ \ie
$\g_{k_\d}< (1/q)^{\a_\d}$.
We now analyze the infinite product (\ref{eq:east4}) which we rewrite as
$$
\prod_{j=k_\d}^\infty\left(\frac{1}{1-\sqrt{\e_j}}\right)\,\prod_{j=k_\d}^\infty\left(1+\frac{1}{s_j}\right).
$$
The second factor, due to the exponential growth of the scales, is
bounded by a constant independent of $q$.
% It is easy to see that
% $$
% j^*=\frac{\log16+\log[1+\frac{\log2}{\log p^{-1}}]}{\log\frac{3}{2}}
% $$

To bound the first factor define
$$
j_*=\min\{j:\e_j\le \nep{-1}\}\approx \log_2(1/q)/(1-\d/2)
$$
and write
\begin{gather}
  \label{eq:east1}
\prod_{j=k_\d}^\infty\left(\frac{1}{1-\sqrt{\e_j}}\right)\le \prod_{j=1}^{j_*}\left(\frac{1+\sqrt{\e_j}}{1-\e_j}\right)
\prod_{j>j_*}^\infty\left(\frac{1}{1-\sqrt{\e_j}}\right)\nonumber\\
\le \nep{C} \,2^{j_*}\,\prod_{j=1}^{j_*}\left(\frac{1}{1-\e_j}\right)
\end{gather}
where we used the bound $1/(1-\sqrt{\e_i})\le
1+\left(e/(e+1)\right)\sqrt{\e_j}$ valid for any $j\ge j_*$
together with
\begin{gather*}
\sum_{j>j_*}^\infty\log\left(1+\frac{\nep{}}{\nep{}+1}\sqrt{\e_j}\right)\le
 \frac{\nep{}}{\nep{}+1}\sum_{j>j_*}^\infty\sqrt{\e_j}\\
\le  \frac{\nep{}}{\nep{}+1} \int_{j_*-1}^\infty dx\
\exp(-q(2^{x(1-\d/2)})/2)
= A_\d\int_{2^{(j_*-1)(1-\d/2)}}^\infty dz\ \exp(-q z/2)/z \\
\le 2A_\d 2^{-(j_*-1)(1-\d/2)}q^{-1}\exp(-q2^{(j_*-1)(1-\d/2)}/2) \le C
\end{gather*}
for some constant $C$ independent of $q$.

Observe now that $1-\e_j\ge 1-\nep{-q\d_j}\ge A q\d_j$ for any $j\le j_*$
and some constant $A\approx \nep{-1}$. Thus the r.h.s. of \eqref{eq:east1}
is bounded from above by
\begin{gather*}
  C\,(\frac{2}{Aq})^{j_*}\,\prod_{j=1}^{j_*}\d_j^{-1}\le
  \frac{1}{q^a}\,(1/q)^{j_*}\,2^{-(1-\d/2)j_*^2/2}
\approx \frac{1}{q^a}\,(1/q)^{\log_2(1/q)/(2-\d)}
\end{gather*}
as $q\downarrow 0$ for some constant $a$.
 \end{proof}

\subsection{FA-1f model}
\label{FA-1f}
In this section we deal with the FA1f model.
Our main results is the following:
\begin{Theorem}
For any $q \in (0,1)$ the spectral gap of the FA-1f model is positive.
\end{Theorem}
\begin{proof}
The proof follows at once from Corollary \ref{IS} because the
probability that the rectangle $\L_0$ of side $\ell$ is internally spanned is
equal to the probability that $\L_0$ is not fully occupied which is equal
to $1-(1-q)^{\ell^d}\uparrow 1$ as $\ell\to \infty$.
\end{proof}
In the next result we discuss the asymptotics of the spectral gap for $
q\downarrow 0$. Such a problem has been discussed at length in the
physical literature with varying results based on numerical simulations
and/or analytical work \cite{WBG1,WBG2,JMS2}.
As a preparation for our bounds we observe that on average
the vacancies are at distance $O(q^{-1/d})$ and each one of them roughly
performs a random walk with jump rate proportional to $q$. Therefore a
possible guess is that
$$
\gap(\cL)= O(q\times \text{gap of a simple RW in a box of side
  $O(q^{-1/d})$ })=O(q^{1+2/d})
$$
Although we are not able to prove or disprove the conjecture for $d\ge
3$ our bounds are consistent with it \footnote{Notice that recent
  work \cite{JMS2} in the physics community
  suggests that $\gap\approx q^2$ for any $d\ge 2$}.
\begin{Theorem}
\label{asisntotica FA-1f}
For any $d\ge 1$, there exists a constant $C=C(d)$
such that for any $q \in (0,1)$, the spectral gap $\gap(\cL)$ satisfies
the following bounds.
$$
\begin{array}{rcccll}
\displaystyle
C^{-1} q^3
& \leq &
\displaystyle \gap(\cL)
& \leq &
C q^3
& \qquad \text{for } d=1, \\
\displaystyle
C^{-1} q^2/\log(1/q) &
\leq &
\displaystyle \gap(\cL)
& \leq &
\displaystyle C q^2
& \qquad \text{for } d=2, \\
\displaystyle
C^{-1} q^2 &
\leq &
\displaystyle \gap(\cL) & \leq
& C q^{1+\frac{2}{d}}
& \qquad \text{for } d \geq 3 .
\end{array}
$$
\end{Theorem}
\begin{proof}
  We begin by proving the upper bounds via a careful choice of a test
  function to plug into the variational characterization for the spectral
  gap.  Fix $d\ge 1$ and assume, without loss
  of generality, $q \ll 1$. Let also $\ell_q =
  \bigl(\frac{\log(1-q_0)}{\log(1-q)}\bigr)^{1/d}\approx
  \l_0 q^{-1/d}$ with $\l_0=|\log(1-q_0)|^{1/d}$, where $q_0$ is as in
  Theorem \ref{main theorem}.

Let $g$ be a smooth function on $[0,1]$ with support in
$[1/4,3/4]$ and such that
\begin{equation}
  \label{eq:tf0}
\int_0^1 \alpha^{d-1} e^{-\alpha^d} g(\alpha) d\alpha = 0
\quad
\text{ and }
\quad
\int_0^1 \alpha^{d-1} e^{-\alpha^d} g^2(\alpha) d\alpha = 1 .
\end{equation}
Set (see figure \eqref{Fig:funztest})
$$
\xi(\sigma) := \sup \left\{ \ell : \sigma(x)=1 \text{ for all }
x \text{ such that } \|x\|_\infty < \ell \right\}
$$
and notice that for any $k=0, \dots, \ell_q$,
\begin{equation}
  \label{eq:tf1}
\mu(\xi = k)=p^{k^d} -p^{(k+1)^d}
\approx
q d k^{d-1} e^{-qk^d}
\end{equation}
Having defined the r.v. $\xi$ the test function we will use is $f=g(\xi/\ell_q)$.
Using \eqref{eq:tf1} together with \eqref{eq:tf0} one can  check
that
\begin{equation}
  \label{eq:tf11}
\Var(f) \approx \frac{1}{\ell_q}\approx q^{1/d}.
\end{equation}
On the other hand, by writing $T_x$ for the spin-flip operator in $x$, {\it i.e.}
\begin{equation*}
  T_x(\s)(y)=
\begin{cases}
\s(y) & \text{if $y\neq x$}\\
1-\s(x) & \text{if $y=x$}
  \end{cases}
\end{equation*}
and using reversibility we have
\begin{eqnarray}
\label{eq:second}
\cD(f)
& = &
\sum_{x\in \bbZ^d} \mu \left[
c_x \left[ g(\frac{\xi \circ T_x}{\ell_q}) - g(\frac{\xi}{\ell_q}) \right]^2
\right]
\nonumber \\
& = &
\sum_{x\in \bbZ^d} \sum_{k=0}^{\ell_q}
\mu \left[
c_x \left[ g(\frac{\xi \circ T_x}{\ell_q}) - g(\frac{\xi}{\ell_q}) \right]^2
\!\!\! \1_{\xi = k} \right] \\
%\mu\left(\{\xi = k \} \right)\\
& = &
2\sum_{k=\inte{\frac14 \ell_q -1}}^{\inte{\frac 34 \ell_q }} %\mu\left(\{\xi = k \} \right)
\Bigl( g\bigl(\frac{k+1}{\ell_q}\bigr) -
g\bigl(\frac{k}{\ell_q}\bigr) \Bigr)^2
\sum_{\substack{x\\\| x \|_\infty=k+1}}
\mu \left( c_x \1_{\xi \circ T_x =k+1} \1_{\xi = k} \right). \nonumber
\end{eqnarray}
Notice that for any $k$, any $x$ such that  $\| x \|_\infty=k+1$,
\begin{gather*}
 \mu \left( c_x \1_{\xi \circ T_x =k+1} \1_{\xi = k} \right)  \\
 \qquad \qquad =
\mu \left( c_x \tc \xi \circ T_x =k+1, \xi = k  \right)
\mu \left( \xi \circ T_x =k+1  \tc \xi = k \right)
\mu\left(\xi = k \right) \\
 \qquad  \qquad \leq
c \, \frac{q}{k^{d-1}} \mu\left(\xi = k  \right)
\end{gather*}
for some constant $c$ depending only on $d$.
The factor $q$ above comes from
the fact that, given $\xi = k$ and $\xi \circ T_x = k+1$,
$x$ is necessarily the only empty site in the $(k+1)^{\rm th}$-layer. Therefore, the flip at
$x$ can occur only if
the nearest neighbor of $x$ in the next layer is empty (see figure
\eqref{Fig:funztest}).
Moreover, given
$\xi = k$, the conditional probability of having
zero at $x$ and the rest of the layer completely filled is
of order $1/k^{d-1}$.
\begin{figure}[t]
\psfrag{k}[][]{{{$k$}}}
\psfrag{l}[][]{{{$k\!\!+\!\!1$}}}
\psfrag{m}[][]{{{$\;\;k\!\!+\!\!2$}}}
\psfrag{x}[][]{{{$x$}}}
\centering
\includegraphics[width=.270\columnwidth]{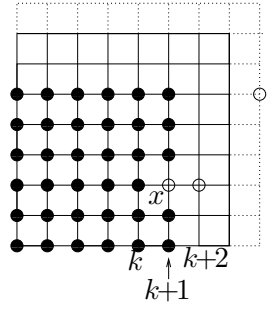}
\caption{In dimension 2, a configuration $\sigma$ where $\xi(\sigma)=k$
and $\xi \circ T_x (\sigma) = k+1$.}
\label{Fig:funztest}
\end{figure}
It follows that
$$
\sum_{x: \| x \|_\infty=k+1}
\mu \left( c_x \1_{\xi \circ T_x =k+1} \1_{\xi = k} \right)
\leq c' q \mu\left( \xi = k \right) .
$$
In conclusion, using \eqref{eq:tf1} and writing $\a:=k/\ell_q$,
\begin{gather}
\cD(f) \leq
c'' q\sum_{k=\inte{\frac14 \ell_q -1}}^{\inte{\frac 34 \ell_q }}
\mu\left( \xi = k \right)\Bigl( g\bigl(\frac{k+1}{\ell_q}\bigr) -
g\bigl(\frac{k}{\ell_q}\bigr) \Bigr)^2 \nonumber\\
\approx
\frac{q}{\ell_q^3} \int_\frac{1}{4}^\frac{3}{4}
\a^{d-1} e^{-(\l_0\a)^d}
g'\left(\a\right)^2
d\a
\approx \frac{q^{1+\frac{2}{d}}}{\ell_q}.
\label{eq:tf4}
\end{gather}
as $q\downarrow 0$.
The upper bound on the spectral gap  follows from
\eqref{eq:tf11},\eqref{eq:tf4} and \eqref{eq:gap}.

\bigskip
We now discuss the lower bound. The first step relates the spectral gap
in infinite volume to the spectral gap in a $q$-dependent finite region.

\begin{Lemma}
  Let $\gap(q)$ be the spectral gap of the FA1f model in
  $\L_{2\ell_q}=\{x\in \bbZ^d: \|x\|_\infty \le 2\ell_q -1\}$ with
  minimal boundary condition, {\it i.e.} exactly one empty site on the
  boundary. There exists a constant $C=C(d)$ such that
$$
\gap(\cL)\ge C \gap(q)
$$
\end{Lemma}
\begin{proof}[Proof of the Lemma]
The starting point is the bound
\eqref{fundam} for $\ell=\ell_q$:
\begin{equation}
    \Var(f)\le 2\sum_{x\in \bbZ^d(\ell_q)}\mu\bigl(\tilde
    c_x\Var_{\L_x}(f)\bigr)
 \label{funda2}
\end{equation}
Recall that $\tilde c_x(\s)$ is simply the indicator
function of the event that for any $y\in \cK^*_{\{x/\ell\}}$ the block
$\L_{\ell y}$
is internally spanned for $\s$ \ie it is not completely filled.
Let us examine a generic term $\mu\bigl(\tilde
c_x\Var_{\L_x}(f)\bigr)$. Given $\s$ such that $\tilde c_x(\s)=1$ let $\xi(\s)$
be the largest $r\le \ell_q$ such that there exists an empty site on
$\partial \L_{x,r}$, where $\L_{x,r}=\{y:\ d_\infty(y,\L_x)\le r\}$.
Exactly as in the proof of Theorem \ref{main theorem} the convexity of
the variance implies that
\begin{equation}
  \mu\bigl(\tilde  c_x\Var_{\L_x}(f)\bigr)\le \mu\bigl(\id_{\xi\le \ell_q}\Var_{\L_{x,\xi}}(f)\bigr)
\label{eq''upb0}
\end{equation}
Since by construction $\Var_{\L_{x,\xi}}(f)$ is computed with an empty site in $\partial
\L_{x,\xi}$, we can use the Poincar\'e inequality for the FA-1f model in
$\L_{x,\xi}$ with \emph{minimal} boundary conditions to get
\begin{equation}
  \mu\bigl(\hat c_x\,\Var_{\L_x}(f)\bigr)\le
  \mu\Bigl(\,\gap(\cL_{\L_{x,\xi}})^{-1}\,
\sum_{z\in \L_{x,\xi}}\mu_{\L_{x,\xi}}\bigl(c_z\Var_z(f)\bigr)\,\Bigr)
\label{eq:upb1}
\end{equation}
By monotonicity of the gap (see Lemma \ref{Mon})
$\gap(\cL_{\L_{x,\xi}})\ge \gap(q)$. Thus the r.h.s. of \eqref{eq:upb1}
is bounded from above by
\begin{equation}
\gap(q)^{-1}\sum_{z\in \L_{x,\ell_q}}\mu\bigl(c_z\Var_z(f)\bigr)
\label{eq:upb1bis}
\end{equation}
If we finally plug \eqref{eq:upb1bis} into  the r.h.s of \eqref{funda2}
we get
$$
\Var(f)\le \gap(q)^{-1} c(d)\sum_{z\in \bbZ^d}\mu\bigl(c_z\Var_z(f)\bigr)
$$
and the Lemma follows.
\end{proof}
The proof of the lower bound will then be complete once we prove the following result.

\begin{Proposition}
There exists a
constant $C=C(d)$ such that
for any $q \in (0,1)$,
\begin{equation}
  \label{eq:upb2}
\gap(q)\ge
C\begin{cases}
  q^3 & \text{if $d=1$}\\
  q^2 /\log(1/q)& \text{if $d=2$}\\
 q^2 & \text{if $d=3$}
\end{cases}
\end{equation}
\end{Proposition}
\begin{proof}[Proof of the proposition]
We begin with the $d=2$ case. For simplicity of notation we simply write
$\L$ for $\L_{2\ell_q}$.

The starting point is the standard Poincar\'e inequality
for the Bernoulli product measure on
$\Lambda$ (see {\it e.g.} \cite[chapter 1]{Bible}).
For  every function $f$
\begin{equation} \label{eq:glauber}
\Var_\L(f) \leq \sum_{x \in \Lambda}
\mu \left( \Var_x (f) \right) .
\end{equation}
Our aim is to bound from above the r.h.s. of \eqref{eq:glauber} with the
Dirichlet form of the FA-1f model in $\L$ with minimal boundary
conditions using a path argument.  Intuitively it works as
follows. Computing the local variance $\Var_x(f)$ at $x$ involves a
spin-flip at site $x$ which might or might not be allowed by the
constraints, depending on the structure of the configuration around
$x$. The idea is then to (see fig. \ref{Fig:path} and \ref{Fig:cammini} for a graphical
illustration):
\begin{enumerate}[(i)]
\item define a geometric path $\g_x$ inside $\L$ connecting $x$ to the (unique) empty
site at the boundary of $\L$;
\item look for the empty site on $\g_x$ closest to $x$;
\item move it,
step by step using allowed flips, to one of the neighbors of $x$ but
keeping the configuration as close as possible to the original one;
\item do the spin-flip at $x$ in the modified configuration.
\end{enumerate}
In order to get an optimal result the choice of the path $\g_x$ is not
irrelevant and we will follow the strategy of \cite{Saloff} to analyze
the simple random walk on the graph consisting of two squares grids sharing exactly one corner.

We first need a bit of extra notation.
We denote by $x^*$ the unique empty site on the boundary $\partial \L$
and for any $y\in \L$ and
any $\h\in \O_\L$ we write $\h^y$ for the flipped configuration
$T_y(\h)$.
Next we declare any pair $e=(\eta,\eta^y)\equiv (e^-,e^+)$ an edge iff
$c_y(\eta)=1$ (\ie the spin-flip at $y$ in the configuration $\h$ is a
legal one). With these notations,
$$
\cD(f) = \sum_e \mu(e^-) \left(f(e^+) - f(e^-) \right)^2 .
$$
To any edge $e=(\eta,\eta^y)$ we associated a
weight $w(e)$ defined by
$w(e)=i+1$ if $d_1(y,x^* )=i$.

Let now, for any
$x \in \Lambda$, $\gamma_x=(x^*,x^{(1)},x^{(2)},\cdots,x^{(n-1)},x)$
be one of the geodesic
paths from $x^*$  to $x$ such that, for any $y \in \gamma_x$, the
Euclidean distance between $y$ and the straight line segment $[x,x^* ]$
is at most $\sqrt 2/2$ (see Figure \ref{Fig:path}).
\begin{figure}[t]
\psfrag{x}[][]{$x$}
\psfrag{o}[][]{$x^*$ }
\centering
\includegraphics[width=.25\columnwidth]{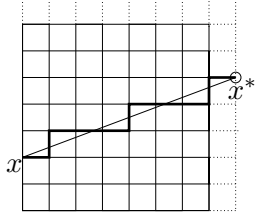}
\caption{An example of geodesic for the path $\gamma_x$.}
\label{Fig:path}
\end{figure}
Given a configuration $\sigma$ we will construct
a path $\Gamma_{\sigma \rightarrow
  \s^x}=\{\h^{(0)},\h^{(1)},\dots,\h^{(j)}\}$, $j\le 2n$, with the properties
that:
\begin{enumerate}[i)]
\item $\h^{(0)}=\s$ and $\h^{(j)}=\s^x$;
\item the path is self-avoiding;
\item for any $i$ the pair $(\h^{(i-1)},\h^{(i)})$ forms an
  edge and the associated spin-flip occurs on $\g_x$;
\item for any $i$ the configuration $\h^{(i)}$ differs from
  $\s$ in at most two sites.
\end{enumerate}
We will denote by
$|\Gamma_{\sigma \rightarrow \sigma^x}|_w:=
\sum_{e \in \Gamma_{\sigma \rightarrow \sigma^x}} \frac{1}{w(e)}$
the weighted length of the path $\Gamma_{\sigma \rightarrow \sigma^x}$.
By the Cauchy-Schwartz inequality, we have
\begin{eqnarray*}
\sum_{x \in \Lambda} \mu \left( \Var_x(f)\right) &= &
pq \sum_{x \in \Lambda} \mu \left( \bigl[f(\s^x)-f(\s)\bigr]^2\right)\\
& =  &
pq\sum_{x \in \Lambda}  \sum_\sigma \mu (\sigma)
\left(
\sum_{e \in \Gamma_{\sigma \rightarrow \sigma^x}}
\frac{\sqrt{w(e)}(f(e^+) - f(e^-))}{\sqrt{w(e)}}
\right)^2 \\
& \leq  & pq \sum_{x \in \Lambda} \sum_\sigma \mu (\sigma)
|\Gamma_{\sigma \rightarrow \sigma^x}|_w \!\!\! \!\!\!\sum_{e \in
\Gamma_{\sigma \rightarrow \sigma^x}}
w(e) \left( f(e^+) - f(e^-) \right)^2 \\
& =  &
pq\sum_{e}
\left( f(e^+) - f(e^-) \right)^2
%\left\{
w(e)
\sum_{\genfrac{}{}{0pt}{}{x \in \Lambda, \sigma :}
{\Gamma_{\sigma \rightarrow \sigma^x} \ni e}}
\mu (\sigma) |\Gamma_{\sigma \rightarrow \sigma^x}|_w
%\right\}
\\
& \leq &
\cD(f) \max_e \Big\{
\frac{pq\,w(e)}{\mu(e^-)}
\sum_{\genfrac{}{}{0pt}{}{x \in \Lambda, \sigma :}
{\Gamma_{\sigma \rightarrow \sigma^x} \ni e}}
\mu (\sigma) |\Gamma_{\sigma \rightarrow \sigma^x}|_w
\Big\} .
\end{eqnarray*}
Fix an edge $e=(\eta,\eta^y)$ with $w(e)=i+1$.
Let $C$ denotes a constant that does not depend on $q$ and
that may change from line to line.
By construction, on one hand we have for any $\sigma$ and $x$ such that
$\Gamma_{\sigma \rightarrow \sigma^x} \ni e$,
$\frac{\mu (\sigma)}{\mu(e^-)} \leq C \frac{1}{q^2}$ because of property
(iii) of $\Gamma_{\sigma \rightarrow \sigma^x}$.
On the other hand, for any  $\sigma$ and $x$,
$$
|\Gamma_{\sigma \rightarrow \sigma^x}|_w
\leq C \sum_{i=1}^{2\ell_q} \frac{1}{i} \leq C \log(\ell_q) .
$$
And finally, by construction, one has (see \cite[section 3.2]{Saloff})
$$
\# \left\{(x,\sigma) : \Gamma_{\sigma \rightarrow \sigma^x} \ni e \right\}
\leq
C \# \left\{ y : \gamma_x \ni y \right\}
\leq
C \frac{|\L|}{i+1} .
$$
Collecting these computations leads to
$$
\sum_{x \in \Lambda} \mu \left(\Var_x(f) \right)
\leq  \frac{C}{q^2} \log(1/q) \, \cD(f) .
$$
\ie the claimed bound on $\gap(q)$.

In $d \geq 3$, the above strategy applies in the same way but one
needs a different choice of
the edge-weight $w(e)$ namely $w(e)=(i+1)^{d-2}$ (see again \cite[section 3.2]{Saloff}).
In $d=1$ instead one can
convince oneself that the  weight function $w \equiv 1$ in the
previous proof leads to the upper bound $1/q^3$, up to some constant.

It remains to discuss the construction of the path $\Gamma_{\sigma
  \rightarrow \sigma^x}$ with the desired properties.  Given $\sigma$,
$x$ and $\gamma_x=\left(x_0=x^*,x^{(1)},x^{(2)},\cdots,x^{(n-1)},x^{(n)}=x\right)$ define
$i_0=\max \{0\leq i \leq n-1 : \sigma(x^{(i)})=0\}$. In this way for any $i
\geq i_0+1$, $\sigma(x^{(i)})=1$. We will denote by $\eta^{x,y}=(\eta^x)^y$
the configuration $\eta$ flipped in $x$ \emph{and} $y$.

If $i_0=n-1$ then trivially
$\Gamma_{\sigma \rightarrow \sigma^x}=\left\{\sigma, \sigma^x\right\}$.
Hence assume that
$i_0 \leq n-2$. We set
$$
\Gamma_{\sigma \rightarrow \sigma^x}=
\left\{\eta^{(0)}=\sigma, \eta^{(1)}, \ldots, \eta^{(2(n-i_0)-1)} =
\sigma^x \right\}
$$
with
$\eta^{(1)}=\sigma^{x^{(i_0+1)}}$ and for
$k=1, \ldots, n-i_0-1$,
$\eta^{(2k)}=\sigma^{x^{(i_0+k)}, x^{(i_0+k+1)}}$,
$\eta^{(2k+1)}=\sigma^{x^{(i_0+k+1)}}$ (see figure \ref{Fig:cammini}).
One can easily convince oneself that $\Gamma_{\sigma \rightarrow
\sigma^x}$
satisfies the prescribed property $(i)-(iv)$ set above.
\begin{figure}[h]
\psfrag{x}{\small$ x$}
\psfrag{y}{\small $x^*$}
\psfrag{i}{\small $x_{i_0}$}
\psfrag{0}{\small $\eta^{(0)}=\sigma$}
\psfrag{1}{\small $\eta^{(1)}$}
\psfrag{2}{\small $\eta^{(2)}$}
\psfrag{3}{\small $\eta^{(3)}$}
\psfrag{4}{\small $\sigma^x$}
\includegraphics[width=.8\columnwidth]{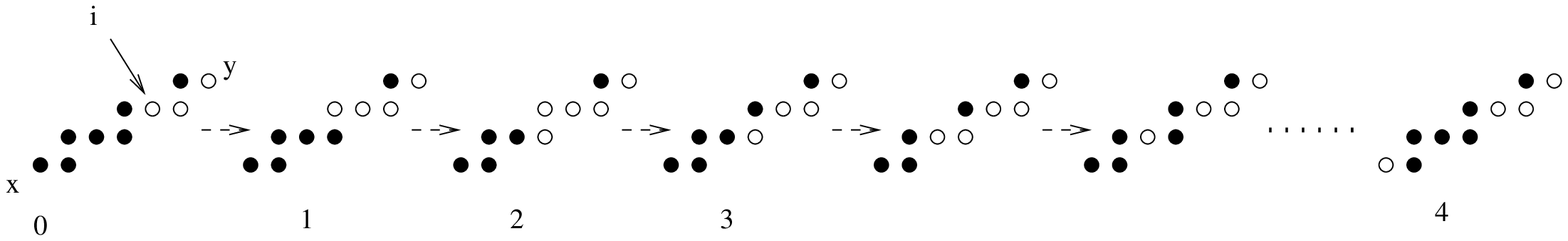}
\caption{The path $\G_{\s\rightarrow \s^x}$.}
\label{Fig:cammini}
\end{figure}
\end{proof}
The proof of the lower bound is complete.
\end{proof}

\subsection{FA-jf and Modified Basic model in $\bbZ^d$}
\label{FA2fmodified}
Next we examine the FA-jf and Modified Basic (MB) model in $\bbZ^d$ with
$d\ge 2$ and $j\le d$.
\begin{Theorem}
For any $q \in (0,1)$ any $d\ge 2$ and $j\le d$ the spectral gap of the
FA-jf and MB models is positive.
\end{Theorem}
\begin{proof}
Under the hypothesis of the theorem both models have a trivial bootstrap
percolation threshold $q_{bp}=0$ and moreover they satisfy the
assumption of corollary \ref{IS} (see \cite{Schonmann}) for any $q>0$. Therefore $\gap
>0$ by Corollary \ref{IS}.
\end{proof}
We now study the asymptotics of the spectral gap as
$q\downarrow 0$ and we restrict ourselves to the most constrained case,
namely either the MB model or the FA-df model. For this purpose we need to introduce few extra
notation and to recall some results from boostrap percolation theory
(see \cite{Holroyd2}).

Let $\d\in \{1,\dots,d\}$. We define the $\d$-dimensional cube
$Q^\d(L):=\{0,\dots,L-1\}^\d\times \{1\}^{d-\d}\sset \bbZ^d$. By a
\emph{copy} of $Q^\d(L)$ we mean an image of $Q^\d(L)$ under any isometry of
$\bbZ^d$.
\begin{definition}
  Given a configuration $\h$, we will say that $Q^\d(L)$ is ``$\d$
  internally spanned'' if $\{1,\dots,L-1\}^\d$ is internally spanned for
  the bootstrap map associated to the corresponding model restricted to
  $\bbZ^{\d}$ (\ie with the rules either of the FA-$\d$f or of the MB
  model in $\bbZ^{\d}$). Similarly for any copy of $Q^\d(L)$.
\end{definition}
Define now
$$
I^d(L,q):=\mu\bigl(\text{$Q^d(L)$ is internally spanned}\bigr)
$$
and let $\exp^n$ denote the $n$-th iterate of the exponential
function. Then the following results is known to hold for both models
\cite{Cerf,Cerf2, Holroyd, Holroyd2}. \\ There exists two positive
constants $0<\l_1\le \l_2$ such that for any $\epsilon >0$
\begin{eqnarray}
\lim_{q\to 0}  I^d\Bigl(\exp^{d-1}(\frac{\l_1-\epsilon}{q}),q\Bigr)&=& 0\\
 \lim_{q\to 0} I^d\Bigl(\exp^{d-1}(\frac{\l_2+\epsilon}{q}),q\Bigr)&=& 1
\label{eq:IS}
\end{eqnarray}
Moreover there exists $c=c(d)<1$ and $C=C(d)<\infty$ such that if $\ell$ is such that
$I^d(\ell,q)\ge c$ then, for any $L\ge \ell$,
\begin{equation}
  \label{eq:IS2}
  I^d(L,q)\ge 1-C\nep{-L/\ell}
\end{equation}
For the FA-2f model and for the MB model for all $d\ge 2$ the threshold
is sharp in the sense that $\l_1=\l_2=\l$ with
$\l=\pi^2/18$ for the FA-2f model and $\l=\pi^2/6$ for the MB model \cite{Holroyd,Holroyd2}.
We are now ready to state our main result.
\begin{Theorem}
Fix $d\ge 2$ and
$\epsilon >0$. Then for both models there exists $c=c(d)$ such that
 \begin{align}
   \label{eq:FA.1}
\left[\exp^{d-1}(c/q^2)\right]^{-1} \le& \gap(\cL) \le
\left[\exp^{d-1}\bigl(\frac{\l_1-\epsilon}{q}\bigr)\right]^{-1} \quad
&d\ge 3\\
\exp(-c/q^5) \le& \gap(\cL) \le
\exp\bigl(-\frac{(\l_1-\epsilon)}{q}\bigr)\quad
&d=2
 \end{align}
as $q\downarrow 0$.
\end{Theorem}
\begin{proof}
  In the course of the proof we will use the following well known
  observation. If a configuration $\h$ is identically equal to $0$ in a
  $d$-dimensional cube $Q$ and each face $F$ of
  $\partial Q$ is ``$(d-1)$ internally spanned'' (by $\h$), then $Q\cup \partial Q$
  is internally spanned.

\medskip
\noindent {\bf (i)}.  We begin by proving the upper bound following the
strategy outlined in remark 3.7. Fix $\epsilon >0$,
let $\L_1$ be the cube centered at the origin of side
$L_1:=\exp^{d-1}\bigl(\frac{\l_1-\epsilon/2}{q}\bigr)$ and let
$m=\exp^{d-2}(\frac{K}{q^2})$ where $K$ is a large constant to be chosen
later on. Define the two events:
  \begin{eqnarray}
A&=&\{\h:\ \L_1 \text{ is not internally spanned}\}\nonumber\\
B &=& \{\h:\ \text{any $(d-1)$-dimensional cube of side $m$ inside
  $\L_1$ is} \nonumber
\\
&\phantom{=}& \text{ ``$(d-1)$ internally spanned''}\}.
\label{eq:B}
  \end{eqnarray}
  Thanks to \eqref{eq:IS} and \eqref{eq:IS2}, $\mu(A)>1/2$ and
  $\mu(B)\ge \frac 34$ if $K$ and $q$ are chosen large enough and small
  enough respectively. Therefore $\mu(A\cap B)\ge 1/4$ for small
  $q$. Pick now $\h\in A\cap B$ and consider $\tilde \h$ which is
  identically equal to one outside $\L_1$ and equal to $\h$ inside. We
  begin by observing that, starting from $\tilde \h$, the little square
  $Q$ of side $m$ centered at origin cannot be completely emptied by the
  bootstrap map $T$ \eqref{eq:bootstrap map}. Assume in fact the
  opposite. Then, after $Q$ has been emptied and using the fact that $\h\in
  B$, we could empty $\partial Q$ and
  continue layer by layer until we have emptied the whole $\L_1$, a
  contradiction with the assumption $\h\in A$.  The above simple
  observation implies in particular that, if we start the Glauber
  dynamics from $\tilde \h$, there exists a point $x\in Q$ such that
  $\s^{\tilde\h}_x(s)=\h_x$ for all $s>0$. However, and this is the
  second main observation, if $t=\frac 14 L_1$, by standard results on
  ``finite speed of propagation of information'' (see
  e.g. \cite{SFlour}) and the basic coupling between the process started
  from $\h$ and the process started from $\tilde \h$,
\begin{gather*}
  \bbP\left(\exists\, x\in Q: \s_x^{\tilde \h}(s)\neq \s_x^\h(s) \text{ for
    some $s\le t$}\right)\ll 1
\end{gather*}
Therefore
\begin{equation*}
  \bbP(\exists\, x\in Q: \; \s_x^\h(s)=\h_x \; \forall s\le \frac 14 L_1)\ge
  \frac 12
\end{equation*}
for all sufficiently small $q$.

We are finally in a position to prove the r.h.s. of
\eqref{eq:FA.1}. Using theorem \ref{persf} combined with the above
discussion we can write
\begin{gather*}
  \nep{-t\frac{q\gap}{2(1+p)}}\ge F(t)\\
\ge \frac{1}{|Q|}\int_{A\cap B} d\mu(\h)\bbP(\exists\, x\in Q: \;
  \s_x^\h(s)=\h_x \; \forall s\le t) \ge \frac {1}{8|Q|}
\end{gather*}
that is $\gap\le c\log\bigl(|Q|\bigr)/qt$ for some constant $c$, \ie the sought
upper bound for $q$ small, given our choice of $t$.

\medskip
\noindent {\bf (ii)} We now turn to the proof of the lower bound in
\eqref{eq:FA.1}. It is enough to consider only the MB model since, being
more restrictive than the FA-df model, it has the smallest spectral
gap.

Fix $\epsilon\in (0,1)$, let
$\ell=\exp^{d-1}\bigl((\l+5\epsilon)/q\bigr)$, $\l=\pi^2/6$, and let $m=
\exp^{d-2}(1/q^2)$ if $d\ge 3$ and $m=K/q^2$ if $d=2$, where $K$ is a
large constant to be fixed later on. Let $E_1$ be the event that
$Q^d(\ell)$ contains some copy of $Q^d(m)$ which is internally spanned
and let $E_2$ be the event that for each $\d\in [1,\dots d-1]$, every copy of
$Q^\d(m)$ in $Q^d(\ell)$ is ``$\d$ internally spanned''. Then it is
possible to show (see section 2 of \cite{Holroyd2} for the case $d\ge 3$
and section 4 of \cite{Holroyd} for the case $d=2$) that both $\mu(E_1)$
and $\mu(E_2)$ tend to one as $q\to 0$ if $K$ is chosen large enough.

Recall now the notation at the beginning of section \ref{main
  results}. The first step is to relate the infinite volume spectral gap
to the spectral gap in the cube $\L_0\equiv Q^d(\ell)$ with zero
boundary condition on $\partial^*_+\L_0$.
  \begin{Proposition}
    There exists a constant $c=c(d)$ such that, for any $q$ small enough,
\[
\gap(\cL)\ge \nep{-cm^d}  \gap(\cL_{\L_0})
\]
  \end{Proposition}
\label{Pro-Fa.1}
  \begin{proof}
    As in the case of the FA-1f model, our starting point is the
    renormalized Poincar\'e inequality \eqref{fundam} on scale $\ell$
    and $\eps_0$-good event $G_\ell:=E_1\cap E_2$. Thanks to
    \eqref{fundam} we can write
\begin{equation*}
  \Var(f)\le 2\sum_{x\in \bbZ(\ell)}\mu\left(\tilde c_x \Var_{\L_x}(f)\right)
\end{equation*}
where the $\tilde c_x$'s are as in \eqref{fundam}.
Without loss of generality we now examine the term $\mu\left(\tilde c_0 \Var_{\L_0}(f)\right)$.
\begin{Lemma}
\label{ultimo}
    There exists a constant $c=c(d)$ such that, for any $q$ small enough,
  \begin{equation*}
    \mu\left(\tilde c_0 \Var_{\L_0}(f)\right)\le
    \nep{cm^d}\gap(\cL_{\L_0})^{-1}\!\!\!\!\!\sum_{x\in \cup_{y\in\cK_0^*
        \cup\{0\}}\L_{\ell y}}\mu\left(c_x\Var_x(f)\right)
  \end{equation*}
where the $c_x$'s are the constraints for the MB
model.
\end{Lemma}
Clearly the Lemma completes the proof of the proposition
  \end{proof}
  \begin{proof}[Proof of the Lemma]
By definition
\[
    \Var_{\L_0}(f)\le \gap(\cL_{\L_0})^{-1}\sum_{x\in
      \L_0}\mu_{\L_0}\bigl(c_{x,\L_0}\Var_x(f)\bigr)
\]
where, we recall, the subscript $\L_0$ in $c_{x,\L_0}$ means that zero
boundary condition on $\partial_+^*\L_0$ are assumed. Notice that, if
$\cK_x^*\sset \L_0$, then $c_{x,\L_0}=c_x$. If we plug the
above bound into $\mu\left(\tilde c_0 \Var_{\L_0}(f)\right)$ and use the
trivial bound $\tilde c_0\le 1$, we see that
all what is left to prove is that
\begin{equation}
  \label{eq:FA.2}
  \mu\bigl(\tilde c_0\, c_{x,\L_0}\Var_x(f)\bigr) \le
\nep{c m^d}\mu\bigl(c_{x}\Var_x(f)\bigr)
\end{equation}
for all $x\in \L_0$ such that $\cK_x\nsubseteq \L_0$.  For simplicity we
assume that $\cK_x\cap \L^c_0$ consists of a unique point $z\in \L_{\ell
  y}$ and we proceed as in the proof of Theorem \ref{main theorem
  01}. Assign some arbitrary order to all cubes of side $m$ inside $
\L_{\ell y}$. Because of the constraint $\tilde c_0$ on the
configuration $\xi$ in $\cup_{y\in \cK_0^*}\L_{\ell y}$, for each $y\in
\cK_0^*$ there exists a sequence of configurations
$(\xi^{(0)},\xi^{(1)},\dots,\xi^{(n)})$, $n\le 2m^d$, with the following
properties:
\begin{enumerate}[(i)]
\item $\xi^{(0)}=\xi$ and $\xi^{(n)}=\xi'$, where $\xi'$ is completely
  empty in the first cube $Q\sset \L_{\ell y}$ of side
  $m$ which was internally spanned for $\xi$ and otherwise coincides
  with $\xi$;
\item $\xi^{(i+1)}$ is obtained from $\xi^{(i)}$ by changing exactly only one
  spin at a suitable site $x^{(i)}\in Q$;
\item the move at $x^{(i)}$ leading from $\xi^{(i)}$ to $\xi^{(i+1)}$ is
  permitted \ie $c_{x^{(i)}}(\xi^{(i)})=1$ for every $i=0,\dots,n$.
\end{enumerate}
\begin{remark}
  Notice that, given $\xi^{(i)}=\h$, the number of starting
  configurations $\xi=\xi^{(0)}$ compatible with $\h$ is bounded from
  above by $2^{cm^d}$, $c=c(d)$, and the relative probability $\mu(\xi)/\mu(\h)$
  by $\bigl(p/q\bigr)^{cm^d}$.
\end{remark}
We can proceed as in \eqref{pippo} and conclude that
\begin{equation}
    \mu\bigl(\tilde c_0\, c_{x,\L_0}\Var_x(f)\bigr) \le
\nep{c' m^d} \mu\bigl(\tilde c_0\,\hat c_0\, c_{x,\L_0}\Var_x(f)\bigr)
\label{eq:FA.3}
\end{equation}
where now $\hat c_0$ is the indicator of the event that for each $y\in
\cK_0^*$ there exists a cube $Q\sset \L_{\ell y}$ of side $m$ which
is completely empty.

Next we observe that for any sequence of adjacent (in e.g. the first
direction) cubes $Q_1,Q_2,\dots Q_j$ of side $m$ inside $\L_{\ell y}$, ordered from
left to right, and for any configuration $\h\in E_2$ which is
identically equal to $0$ in $Q_1$, one can construct a sequence
of configurations $(\h^{(0)},\h^{(1)},\dots,\h^{(n)})$, $n\le j m^d$,
such that:
\begin{enumerate}[(i)]
\item $\h^{(0)}=\h$ and $\h^{(n)}$ is completely
  empty in $Q_j$  and otherwise coincides with $\h$;
\item $\h^{(i+1)}$ is obtained from $\h^{(i)}$ by changing exactly only one
  spin at a suitable site $x^{(i)}\in \cup_{i=1}^jQ_i$;
\item the move at $x^{(i)}$ leading from $\h^{(i)}$ to $\h^{(i+1)}$ is
  permitted \ie $c_{x^{(i)}}(\h^{(i)})=1$ for every $i=0,\dots,n$.
\end{enumerate}
In other words one can move the empty square $Q_1$ to the position
occupied by $Q_j$ in no more than $jm^d$ steps.  The construction is
very simple and it is based on the basic observation described at the
beginning of the proof. Starting from $Q_1$ and using the fact that any
copy of $Q^{d-1}(m)$ inside $\L_{\ell y}$ is ``$(d-1)$ internally
spanned'', by a sequence of legal moves one can first empty $Q_2$. Next one
repeats the same scheme for $Q_3$. Once that also $Q_3$ has been
emptied one backtracks and readjust all the spins inside $Q_2$ to their
original value in the starting configuration $\h$. The whole procedure
is then iterated until the last square $Q_j$ is emptied and the
configuration $\h$ fully reconstructed in $\cup_{i=1}^{j-1}Q_i$.

The key observation at this point is that, given an intermediate
step $\h^{(i)}$ in the sequence, the number of starting
configurations $\h$ compatible with $\\h^{(i)}$ is bounded from
above by $2j\cdot 4^{m^d}$ and the relative probability
$\frac{\mu(\h^{(i)})}{\mu(\h)}$ by $\bigl(p/q\bigr)^{2m^d}$.

By using the path argument above and by proceeding again as in \eqref{pippo},
we can finally bound from above the r.h.s. \eqref{eq:FA.3} by
\begin{equation*}
2\ell\, \nep{c'' m^d}
\mu\bigl(\tilde c_0\,\hat c_{0,x}\, c_{x,\L_0}\Var_x(f)\bigr)
\end{equation*}
where $\hat c_{0,x}$ is the indicator of the event that there exists a
cube $Q$ of side $m$, laying outside $\L_0$ but such that
$\cK_x^*\cap \L_0^c\sset Q$ , which is
completely empty.  Clearly $\tilde c_0\,\hat c_{0,x}\, c_{x,\L_0}\le
c_x$ because the sites in $\cK_x^*\cap \L_0^c$ are forced to be empty
and the proof of the Lemma is complete.
  \end{proof}
  As a second step we lower bound $\gap(\cL_{\L_0})$ by the spectral gap
  in the reduced volume $\L_1:=Q^d_{\ell/2}$ (we assume here for
  simplicity that both $\ell$ and $m$ are powers of $2$). To this end we
  partition $\L_0$ into disjoint copies of $\L_1$,
  $\{\L_1^{(i)}\}_{i=1}^{2^d}$ and, mimicking the argument of section
  \ref{GM}, we run the constrained dynamics of the $*$-general model on
  $\L_0$ with blocks $\{\L_1^{(i)}\}_{i=1}^{2^d}$ and good event the
  event that for each $\d\in [1,\dots d-1]$, every copy of $Q^\d(m)$ in
  $Q^d(\ell/2)$ is ``$\d$ internally spanned''.  By choosing the
  constant $K$ appearing in the definition of $m$ larfge enough the
  probability of $G$ is very close to one as $q\to 0$ and therefore the
  Poincar\'e inequality
  \begin{equation}
    \label{eq:FA8}
    \Var_{\L_0}(f)\le 2\sum _{i=1}^{2^d}\mu\bigl(c_i\Var_{\L_1^{(i)}}(f)\bigr)
  \end{equation}
holds, where $c_i$ are the constraints of the $*$-general model. At this
point we can proceed exactly as in the proof of lemma \ref{ultimo}
and get that the r.h.s. of \eqref{eq:FA8} is bounded from above by
\begin{equation*}
  \nep{c m^d}\gap(\cL_{\L_1})^{-1}\cD_{\L_0}(f)
\end{equation*}
for some constant $c=c(d)$. We have thus proved that
$$
\gap(\cL_{\L_0})^{-1}\le  \nep{c m^d}\gap(\cL_{\L_1})^{-1}
$$
If we iterate $N$ times, where $N$ is such that $2^{-N}\ell=m$ we
finally get
$$
\gap(\cL_{\L_0})^{-1}\le  \nep{c N\,m^d}\gap(\cL_{\L_N})^{-1}
$$
where $\L_N=Q^d_{m}$.
\end{proof}

\subsection{The N-E model}
The N-E model is the natural two dimensional analogue of the
one dimensional East model.
Before giving our results we need to
recall some definitions of the oriented
percolation \cite{Durrett,Schonmann}.
A \emph{NE oriented path} is a collection
$\{x^{(0)},x^{(1)},\cdots,x^{(n)}\}$ of distinct points in $\bbZ^2$
such that
$x^{(i+1)}=x^{(i)}+\a_1\vec e_1+\a_2\vec e_2,~\a_j=0,1~ {\rm and}~\a_1+\a_2=1$
for all $i$. Given a configuration $\eta\in \O$ and $x,y\in\bbZ^2$,
we say that $x\rightarrow y$ if there is a NE
oriented path of occupied sites starting in $x$ and ending in $y$.
For each site $x\in\bbZ^2$ its
\emph{NE occupied cluster} $x$ is the random set
\begin{equation*}
C_{x}(\eta):=\{y\in\bbZ^2~:~~x\rightarrow y\}
\end{equation*}
The range of $C_x(\eta)$ is the random variable
\begin{equation*}
  A_x(\h)=
  \begin{cases}
    0 & \text{if $C_x(\eta)=\emptyset$}\\
 \sup\{1+\|y-x\|_1\,:\,y\in C_x(\eta)\} & \text{otherwise}
  \end{cases}
\end{equation*}
\begin{remark}
If $A_x(\h)>0$ then at least $A_x(\h)$ legal (\ie fulfilling the NE
constraint) spin flip moves are needed to empty the site $x$.
\end{remark}
Finally we define the monotonic
non decreasing function $\theta(p):=\mu(A_0=\infty)$ and let
$$
  p_c^o=\inf\{p\in[0,1]\,:\,\theta(p)>0\}
$$
It is known (see \cite{Durrett}) that $0<p_c^o<1$.  In \cite{Schonmann}
it is proven that the percolation threshold and bootstrap percolation
threshold (see section 2.3) are related by $p_c^o=1-q_{bp}$ and
therefore, thanks to proposition \ref{qc=qbp}, $q_c=1-p_c^o$.  The
presence of a positive threshold $q_c$ reflects a drastic change in the
behavior of the NE process when $q<q_c$ due to the presence of blocked
configurations (NE occupied infinite paths) with probability one.  In
\cite{Lalley} it is proven that the measure $\mu$ on the configuration
space is mixing for $q\ge q_c$, a result that follows at once from the
the arguments given in the proof of proposition \ref{qc=qbp} since
$\theta(p_c^o)=0$ \cite{Grimmett2}.

We now analyze the spectral gap of the N-E process above, below and at
the critical point $q_c$.

\\
Case ${\bf q>q_c}$. This region is characterized by the following result
of \cite{Durrett}.
\begin{Proposition}\label{propop1}
  If $p<p_c^o$ there exists a positive constant
  $\varsigma=\varsigma(p)>0$ such that
\begin{equation}\label{expdecay}
\lim_{n\rightarrow\infty}-\frac{1}{n}\log\mu(A_0\ge\,n)=\varsigma
\end{equation}
\end{Proposition}
We can now state our main theorem
\begin{Theorem}\label{ne1}
For any $q>q_c$ the spectral gap of N-E model is positive.
\end{Theorem}
\begin{proof}
Recall the notation of section \ref{main results}.
Using theorem \ref{main theorem 01} we need to find a set of configurations
$G_\ell$ satisfying properties (a) and (b) of definition \ref{defgelle}.
Fix $\d\in(0,1)$ and $\ell>2$ and define
$$
G_{\ell}:=\{\eta\in\{0,1\}^{\L_0}:\nexists\text{ occupied oriented path in
  $\L_0$ longer than } \ell^\d\}
$$
Since $q>q_c$ we can use (\ref{expdecay}) to obtain that for any
$\epsilon\in (0,1)$ there exists $\ell_c(q,\e,\d) $ such that, for any $
\ell\ge\ell_c(q,\e,\d)$, $\mu(G_\ell)\ge
1-\e$ and property (a) follows. Property (b) also follows directly from the
definition of $G_\ell$. Indeed, if the restriction of a configuration $\h$ to
each one of the squares $\L_0+\ell x$, $x\in \cK_0^*$, belongs to
$G_\ell$, then necessarily there is no occupied oriented path in
$\cup_{x\in \cK_0^*}\{\L_0+\ell x\}$ of length greater than $3\ell^\d$. Therefore, by a sequence of legal
moves, all the $\partial_+^*\L_0$ can be emptied for $\h$ and the proof
is complete.
\end{proof}\noindent
Case ${\bf q<q_c}$.
Following \cite{Durrett} we need few extra notation.
For every $L\in\bbN$ and $\h\in \O$ let $C_0^{(L)}(\h)=\{x\in C_0(\h):
\|x\|_1=L\}$ and let
$$
\xi_0^{(L)}(\h):=\cup_{x\in C_0^{(L)}(\h)}\{x_1\}
$$
be the projection onto the first coordinate axis of
$C_0^{(L)}(\h)$. Denote by $r_L,l_L$ the
right and left edge of $\xi_0^{(L)}(\h)$ respectively.
If $p>p_c^o$ it is possible to show \cite{Durrett} that there exists
positive constants $a, \zeta$ such that
\begin{equation}
  \label{eq:lr}
\mu\left(\{\xi_0^{(L)}\neq \emptyset\} \cap \{r_L\le aL\}\right)=
\mu\left(\{\xi_0^{(L)}\neq \emptyset \}\cap \{l_L\ge aL\}\right)\le \nep{-\zeta L}
\end{equation}
for any $L$ large enough.
\begin{figure}[h]
\psfrag{0}{$0$}
\psfrag{L}{$L$}
\psfrag{x1}{$x_1$}
\psfrag{x2}{$x_2$}
\psfrag{lL}{$l_L$}
\psfrag{rL}{$r_L$}
\psfrag{s}{$\left\{x: \vert| x \vert_1= L \right\}$}
\includegraphics[width=.80\columnwidth]{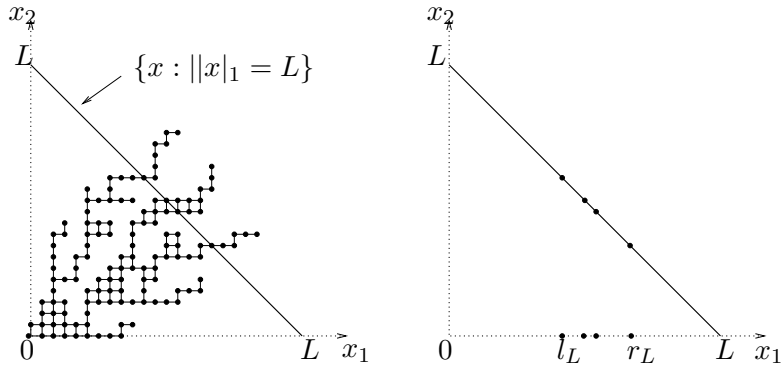}
\caption{
An example of configuration $\eta$ with the sets
$C_0(\eta)$ (on the left) , $C_0^{(L)}(\eta)$ and
$\xi_0^{(L)}(\eta)$ (on the right).}
\end{figure}
We can now state our result for the spectral gap.
\begin{Theorem}\label{ne2}
Let $\L\subset\bbZ^2$ be a square of side $L\in\bbN$.
For any $q<q_c$ there exists two positive constants $c_1$,
$c_2$ such that
\begin{equation}
  \label{eq:ul}
\exp\{-c_1L\}\le{\rm gap}({\cL}_\L)\le\exp\{-c_2L\}
\end{equation}
\end{Theorem}
\begin{proof}
  We first discuss the upper bound by exhibiting a suitable test
  function $f$ to be plugged into the variational characterization of the
  spectral gap.  For this purpose let $
  B_L:=\{\eta:\xi_0^{(L)}\neq\emptyset\} $ and define
  $f=\id_{B_L}$. Since $q<q_c$, there
  exists two positive constants $0<k_1(q)\le k_2(q)<1$ such that
  $k_1\le\mu(B_L)\le k_2$ , see
  \cite{Durrett}.  Thus the variance
  of $f$ is bounded from below uniformly in $L$. On the other hand, by construction,
$$
\cD(f)=\sum_{x\in\L}\mu\left(c_x\Var_x(f)\right)\le\,
|\L|\,\mu({\bar B}_L)
$$
where ${\bar B}_L:=\{\eta\,:\,|\xi^{(L)}_0|=1\} = \{\eta\,:\,r_L=l_L\}$.
Thanks to \eqref{eq:lr}
\begin{eqnarray*}
\mu({\bar B}_L)&
\le\,\mu(\{r_L=l_L\}\cap\{r_L>aL\})+\mu(\{r_L=l_L\}\cap\{r_L\le aL\})\\
&\le\,2\mu(\xi_0^{(L)}\neq \emptyset\}\cap \{r_L\le aL\})\le\,\exp\{-\zeta
L\}
\end{eqnarray*}
and the r.h.s. of \eqref{eq:ul} follows.

\noindent
The bound from below comes from the bisection method of theorem
\ref{main theorem}
where in proposition \ref{gapblock} $\e_k$ is defined as
 the probability
that there is at least one left-right NE occupied oriented path.
Trivially $\e_k\le 1-e^{-c\d_k}$ for some constant $c$. If we plug such
a bound
bound into (\ref{eq:4}) and we remember that the number of steps
of the iterations grows as $c\log L$, we obtain the desired result.
\end{proof}
\noindent
The case ${\bf q=q_c}$.
\begin{Theorem}
The spectral gap is continuous at $q_c$ where, necessarily, it is
zero.
\end{Theorem}
\begin{proof}
Assume $q=q_c$ and suppose that the spectral gap is positive.
Then, by Theorem \ref{persf}, the
persistence function decays exponentially fast as $t\to \infty$. We will show that such a
decay necessarily implies that the all moments of the size of the
oriented cluster $C_0$ are finite \ie $q>q_c$, a contradiction.

Let $H(t):=\{\eta\,:\, A_0(\eta)\ge\,2t\}$ and observe that, again by
the ``finite speed of propagation'' (see section \ref{FA2fmodified}),
$\bbP(\s^\h_0(s)=\h_0\ \text{for all $s\le t$})\ge \frac{1}{2}$ for all
$\h\in H(t)$. Using $H(t)$ we can lower bound $F(t)$ as follows.
\begin{eqnarray*}
F(t)&=&\int\, d\mu(\h)\, \bbP(\s^\h_0(s)=\h_0\ \text{for all $s\le t$})\\
&\ge&\,
\int_{H(t)}\, d\mu(\h)\,\bbP(\s^\h_0(s)=\h_0\ \text{for all $s\le t$})\\
&\ge&\frac{1}{2}\mu(A_0\ge\,2t)
\end{eqnarray*}
which implies,
\begin{equation}\label{eqnpers}
\mu(A_0\ge\,2t)\le\,2 F(t)\le\,2\,e^{-ct}
\end{equation}
for a suitable constant $c>0$.
But \eqref{eqnpers} together with the fact that $|C_0|\le A_0^2+1$ implies that $\mu\left(|C_0|^n\right)<\infty$ for
all $n\in \bbN$, \ie
$p<p_c^o$ \cite{AB}.

The same argument proves continuity at $q_c$.\\
 Suppose in fact that
$\limsup_{q\downarrow q_c} \gap >0$. That would imply \eqref{eqnpers} for any
$q>q_c$ with $c$ independent of $q$, \ie
$\sup_{q>q_c}\mu\left(|C_0|\right)<\infty$, again a contradiction since
$\mu\left(|C_0|\right)$ is an increasing function of $q$ which is
infinite at $q_c$ \cite{Durrett, Grimmett}.
\end{proof}
% \begin{remark} It is usually conjectured that the correlation length $\varsigma$
% diverges at the percolation threshold with a critical exponent, i.e.
% there exists $\rho>0$ such that $\varsigma\approx |p-p_c^o|^{-\rho}$.
% If this is true from (\ref{expdecay}) and
% (\ref{eqnpers}) it follows that near the threshold the
% spectral gap is such that $\gamma(q)\le |q-q_c|^\rho$.
% \end{remark}

\begin{Corollary}
\label{slow}
At $q=q_c$ the persistence function $F$ satisfies
$$
\int_0^\infty dt \, F(\sqrt{t})=\infty
$$
\end{Corollary}
\begin{proof}
By \eqref{eqnpers}
\begin{gather*}
\int_0^\infty dt\, F(\sqrt{t})\ge
\frac 12 \int_0^\infty\!\! dt\, \mu(A_0\ge\,\sqrt{2t})\\
\ge \frac 12  \int_0^\infty dt\, \mu\left(|C_0|\ge c' t\right) =+\infty
\end{gather*}
because  $\mu\left(|C_0|\right)=+\infty$ at $q_c$.
\end{proof}
%\subsection{The Knights}

\section{Some further observations}
\label{Some further observations} We collect here some further comments and aside results that so far have
been omitted for clarity of the exposition.

\subsection{Logarithmic and modified-logarithmic Sobolev constants}
\label{entropy}\\
 A first natural question is whether it would be
possible to go beyond the Poincar\'e inequality and prove a stronger
coercive inequality for the generator $\cL$ like the
\emph{logarithmic} or \emph{modified-logarithmic} Sobolev
inequalities \cite{Bible}. As it is well known, the latter is weaker
than the first one and it implies in particular that, for any
non-negative mean one function $f$ depending on finitely many
variables, the entropy $\ent(P_tf):=\mu\left(P_tf\log(P_tf)\right)$
satisfies:
\begin{equation}
  \label{eq:entropy}
  \ent(P_tf)\le \ent(f)\nep{-\a t}
\end{equation}
for some positive $\a$.  As we briefly discuss below such a behavior is
in general impossible and both the (infinite volume) logarithmic and
modified logarithmic Sobolev constants are zero\footnote{In finite
  volume with minimal boundary conditions it is not difficult to show
  that for some of the models discussed before the logarithmic Sobolev
  constant shrinks to zero as the inverse of the volume}. For simplicity
consider any of the \text{0-1 KCSM} analyzed in section
\ref{specificmodels} and choose $f$ as the indicator function of the
event that the box of side $n$ centered at the origin is fully occupied,
normalized in such a way that $\mu(f)=1$. Denote by $\mu^f$ the
probability measure whose relative density w.r.t. $\mu$ is $f$. If we
assume \eqref{eq:entropy} the relative entropy $\ent(\mu^fP_t/\mu)$
satisfies
\begin{equation}
  \label{eq:entropy2}
\ent(\mu^fP_t/\mu)=\ent(P_tf)\le Cn^d\nep{-\a t}.
\end{equation}
which implies, thanks to Pinsker inequality, that
\begin{equation}
  \label{eq:entropy3}
\|\mu^fP_t -\mu P_t\|^2_{TV}=\|\mu^fP_t -\mu\|^2_{TV}\le
2\ent(\mu^fP_t/\mu)\le 2Cn^d\nep{-\a t}
\end{equation}
\ie $\|\mu^fP_t -\mu P_t\|_{TV}\le \nep{-1}$ for any $t\ge
O(\a^{-1}\log(n))$. However the above conclusion clashes with a standard
property of interacting particles systems with bounded rates known as
``finite speed of propagation'' (see e.g. \cite{SFlour}) which can be
formulated as follows. Let $\t(\h)$
be the first time the origin is updated starting from the configuration
$\h$. Then $\int d\mu^f(\h)\bbP(\t(\h) < t)\le C n^{d-1}\bbP(Z \ge n/r)$
where $Z$ is a Poisson variable of mean $t$ and $r$ is the range defined
in section 2.2.  The above bound implies in particular that $\int
d\mu^f(\h)\bbE(\s^\h_0(t))\approx 1$ for any $t\ll n$ \ie a
contradiction with the previous reasoning.

\subsection{More on the ergodicity/non ergodicity issue in finite
volume}\\
In section \ref{setting} we mentioned that one could try to analyze
a 0-1 KCSM in a finite region without inserting appropriate boundary
conditions in order to guarantee ergodicity but rather by
restricting the configuration space to a suitable ergodic component.
Although such an approach appears rather complicate for e.g.
cooperative models, it is within reach for non-coop\-er\-ative
models.

For simplicity consider the FA-1f model in a finite interval
$\L=[1,\dots,L]$ with configuration space $\O^+_\L :=\{\h\in
\O_\L:\ \sum_{x\in \L}\h_x< L\}$, \ie configuration with at least one
empty site, and constraints corresponding to boundary conditions outside
$\L$ identically equal to one. In other words the constraints only
consider sites inside $\L$. The
resulting Markov process is ergodic and reversible w.r.t the
conditional measure $\mu_\L^+:=\mu_\L(\cdot\tc \O^+_\L)$. We now show how to
derive from our previous results that also the spectral gap of this new
process stays uniformly positive as $L\to \infty$. To keep the
notation simple we drop the subscript $\L$ from now on.

\smallno
For any $\h\in \O^+$, let $\xi(\h)=\min\{x\in \L:\h_x=0\}$ and write,
for an arbitrary $f$,
\begin{equation}
  \label{eq:more}
  \Var^+(f)=\mu^+\left(\Var^+(f\tc \xi)\right)
  +\Var^+\left(\mu^+(f\tc \xi)\right)
\end{equation}
with self explanatory notation. Since $\Var^+(f\tc \xi)$ is computed
with ``good'', \ie zero, boundary condition at $\xi$, we  get that
\begin{gather}
  \label{eq:more2}
  \Var^+\left(f\tc \xi)\right)=\Var\left(f\tc \xi)\right)\\ \le
  {\rm const.}
  \sum_{x<\xi}\mu\left(c_x \Var_x(f) \tc
    \xi\right)={\rm const.}\sum_{x<\xi}\mu^+\left(c_x \Var_x(f) \tc \xi\right).
\end{gather}
Therefore the first term in the r.h.s of \eqref{eq:more} is bounded from
above by a constant times the Dirichlet form. In order to bound the second term in
the r.h.s of \eqref{eq:more} we observe that $\xi$ is a geometric random
variable condition to be less or equal than $L$. By the classical
Poincar\'e inequality for the geometric distribution, we can then write
\begin{gather}
  \Var^+\left(\mu^+(f\tc \xi)\right)\nonumber \\\le {\rm const.}\sum_{x=1}^{L-1}
\mu^+\left(b(x)\left[\mu^+(f\tc \xi=x+1)-\mu^+(f\tc \xi=x)\right]^2\right)
\label{more3}
\end{gather}
where $b(x)=\mu^+(\xi=x+1)/\mu^+(\xi=x)$.
A little bit of algebra now shows that
\begin{gather}
  \mu^+(f\tc \xi=x)-\mu^+(f\tc
  \xi=x+1)=\nonumber \\ =\mu^+\left(\h_{x+1}(f(\h)-f(\h^{x+1})\tc \xi=x\right)
+ \mu^+\left(f(\h^x)-f(\h)\tc \xi=x+1\right)\nonumber\\
=\mu^+\left(\h_{x+1}c_{x+1}(f(\h)-f(\h^{x+1})\tc \xi=x\right)\nonumber\\
+ \mu^+\left(c_x(f(\h^x)-f(\h))\tc \xi=x+1\right)
\label{more4}
\end{gather}
In the last equality we have inserted the constraints $c_{x+1}$ and
$c_x$ because they are identically equal to one. If we now insert
\eqref{more4} into the r.h.s. of \eqref{more3} and use Schwartz
inequality, we get that also the second term in the r.h.s of
\eqref{eq:more} is bounded from above by a constant times the Dirichlet
form and the spectral gap stay bounded away from zero uniformly in $L$.

%\section{Conclusions and open problems}

\bibliographystyle{amsplain}
\bibliography{ref}

\\

\end{document}